\documentclass[,nonblindrev,11pt]{IIETransactions} 
\usepackage[margin=1in]{geometry}
\usepackage{algorithm,standalone}
\usepackage{algpseudocode}
\usepackage{wrapfig,enumerate,enumitem, multirow}
\usepackage{tabularx} 
\usepackage{tabulary,booktabs, threeparttable, stackengine,array}
\usepackage{natbib,textcomp}
\usepackage{graphicx}
\usepackage[toc,page,title,titletoc]{appendix}
\usepackage{tikz}
 \usepackage[scr=euler]{mathalfa}
 \usepackage{bbm}
 
 \usepackage{amsmath}
 \usepackage{amsfonts}
 \usepackage{caption}
 \usepackage{graphicx}
 \usepackage{subcaption}
 \usepackage{algorithm, algpseudocode}
 \usepackage{pgfplots}
 \pgfplotsset{compat=1.8}
 \usepgfplotslibrary{statistics}
 \usepackage{pgf} 
 \usepackage{pgfplots}
 \usetikzlibrary{fit,calc}
 \usepgfplotslibrary{external}
 \usetikzlibrary{decorations.pathmorphing}
 \usetikzlibrary{patterns}
 \newcommand{\boxplot}[6]{%
 	\filldraw[fill=white,line width=0.2] let \n{boxxl}={#1-0.1}, \n{boxxr}={#1+0.1} in (axis cs:\n{boxxl},#3) rectangle (axis cs:\n{boxxr},#4);   
 	\draw[line width=0.2mm, color=red] let \n{boxxl}={#1-0.1}, \n{boxxr}={#1+0.1} in (axis cs:\n{boxxl},#2) -- (axis cs:\n{boxxr},#2);             
 	\draw[line width=0.2mm] (axis cs:#1,#4) -- (axis cs:#1,#6);                                                                           
 	\draw[line width=0.2mm] let \n{whiskerl}={#1-0.025}, \n{whiskerr}={#1+0.025} in (axis cs:\n{whiskerl},#6) -- (axis cs:\n{whiskerr},#6);        
 	\draw[line width=0.2mm] (axis cs:#1,#3) -- (axis cs:#1,#5);                                                                           
 	\draw[line width=0.2mm] let \n{whiskerl}={#1-0.025}, \n{whiskerr}={#1+0.025} in (axis cs:\n{whiskerl},#5) -- (axis cs:\n{whiskerr},#5);        
 }
 \usepackage{rotating}

\usetikzlibrary{shapes.geometric, arrows}
\definecolor{mygreen}{RGB}{0,100,0}
\tikzset{snake it/.style={decorate, decoration={snake, amplitude=0.3mm, segment length=2mm}}}
\tikzstyle{process} = [rectangle, minimum width=3cm, minimum height=1cm, text centered, draw=black]
\tikzstyle{decision} = [diamond, minimum width=3cm, minimum height=1cm, text centered, draw=black,aspect=4]
\tikzstyle{startstop} = [rectangle, rounded corners, minimum width=2cm, minimum height=1cm,text centered, draw=black]
\tikzstyle{io} = [trapezium, trapezium left angle=70, trapezium right angle=110, minimum width=3cm, minimum height=1cm, text centered, draw=black]
\tikzstyle{arrow} = [thick,->,>=stealth]
\tikzstyle{line} = [draw, thick, -latex']

\newcommand*\xbar[1]{%
  \hbox{%
    \vbox{%
      \hrule height 0.5pt 
      \kern0.5ex
      \hbox{%
        \kern-0.2em
        \ensuremath{#1}%
        \kern-0.1em
      }%
    }%
  }%
}

\newcolumntype{C}{>{\centering\arraybackslash}m{3cm}}
\newcolumntype{L}{>{\arraybackslash}m{2.5cm}}

 \bibpunct[, ]{(}{)}{,}{a}{}{,}%
 \def\bibfont{\small}%
 \def\bibsep{\smallskipamount}%
 \def\bibhang{24pt}%
\usepackage{float}
\TheoremsNumberedThrough     
\EquationsNumberedThrough    
\ARTICLEAUTHORS{%
\AUTHOR{Sandra D. Eksioglu$^1$, Berkay Gulcan$^2$, Mohammad Roni$^3$, Scott Mason$^2$}
\AFF{$^1$University of Arkansas, \EMAIL{sandra@uark.edu}}
\AFF{$^2$Clemson University, $^3$Idaho National Laboratory}
} 
\RUNAUTHOR{Eksioglu et al.}
\RUNTITLE{Biomass Blending}

 \usepackage{verbatim}

\begin{document}

\TITLE{\large A Stochastic  Biomass Blending Problem in Decentralized Supply Chains}

\ABSTRACT{%
{Blending biomass materials of different physical or chemical properties provides an opportunity to  adjust the quality of the feedstock to meet the specifications of the conversion platform. We propose a model which identifies the right mix of biomass to optimize the performance of the Thermochemical conversion process at the minimum cost. This is a chance-constraint programming (CCP) model which takes into account the stochastic nature of biomass quality. The proposed CCP model ensures that process requirements, which are impacted by  physical and chemical properties of biomass, are met most of the time. We consider two problem settings, a centralized and a decentralized supply chain. We propose a  mixed-integer linear program to model the blending problem in the centralized setting and a bilevel program to model the blending problem in the decentralized setting.  We use the sample average approximation (SAA) method to approximate the chance constraints, and propose solution algorithms to solve this approximation. We develop a case study for South Carolina  using data provided by the Billion Ton Study. Based on our results, the blends identified consist mainly of pine and softwood residues. The cost of the centralized supply chain is 2 to 6\% lower, which shows that the assumption of centralized decision making leads to  underestimating costs in the supply chain.}
%
%
}%
\maketitle

\section{Introduction}
The majority of the existing biomass logistics models are focused on reducing the total costs of delivering biomass to conversion plants. This is mainly because the emerging biomass supply system inherited models (and the underlying assumptions) from the existing agricultural and logging industries. These models pay little attention on the impacts that biomass quality has on costs and conversion rates \citep{kenney2013understanding}.

Typically, feedstocks of high quality are expensive. Blending biomass feedstock of different physical or chemical properties provides an opportunity to  adjust the quality of the feedstock to meet annual needs of the conversion platform at the minimum cost.  Quality indicators for biomass are moisture content, thermal content, ash content, etc.  For example, clean pine is frequently identified as the biomass feedstock of choice for Thermochemical conversion process. Factors which impact the performance of pyrolysis are  oxygen and hydrogen content of biomass. Blending clean pine in appropriate proportions with logging residues would result in a blend that maintains the desired levels of ash content, while using relatively less expensive materials. Blending for such purposes is a common practice in many industries \citep{hill1990grain}.
Other examples are blending of animal feed to obtain specific nutrient requirements \citep{reddy2009precision}; and blending of high-ash biomass with low-ash coal to allow their use for biopower generation \citep{sami2001co}.


We propose a stochastic optimization model that identifies a blending of different  types of biomass to meet process requirements of a biorefinery at the minimum cost. Current processes require that  ash content of the blend used must be lower than a threshold value, and  total thermal  content must be higher than another threshold value. Meeting these requirements is a challenge because ash and thermal contents of biomass are random and vary by supplier. These process requirements are soft, since, for example, if ash content is higher than the threshold value, the blend undergoes preprocessing  which reduces ash content. Similarly, if thermal content is low, then, additional biomass can be purchased to meet the thermal requirement by contracting new suppliers. However, these practices are expensive. Thus, ash and thermal requirements  should be met most of the time (e.g., 80-90\% of the time) because meeting these requirements all the time impose high costs to the system. We model these soft requirements using chance constraints. 

US Department of Energy and other federal agencies have made significant investments to help the bioenergy industry grow. Despite of these investments, this industry remains a nascent concern and unable to compete with fossil fuels. The purpose of this work is to evaluate the potential impact that practices, such as biomass blending, have to reduce the cost of producing biofuels. To achieve this, we propose two supply chain models, one that assumes centralized decision making in which a single decision maker has full control, and another model that assumes decentralized decision making in a noncooperative environment. The decentralized model considers that the biorefinery and suppliers are independent entities  who have their own goals and objectives. In this model, each entity makes decisions to improve its own performance rather than the performance of the overall supply chain.  To address this, we model the relationships between the biorefinery and suppliers via a Stackelberg game.  The goal is to observe the impact that the assumption of centralized decision making has on the blends identified and corresponding costs.

We develop a case study focused on South Carolina (SC). To develop this case study we used the Billion Ton Study \citep{langholtz2016} which provides data about the availability and cost of different types of biomass in each county of SC. This cost includes the cost of land, and the cost if planting, harvesting, and collection biomass. We also used the Bioenergy Feedstock Library \citep{BFL} which is developed and maintained by the Idaho National Laboratory (INL) with sponsorship from the U.S. Department of Energy (DOE). This library is both a physical repository and knowledge database of biomass feedstock and provided the data about ash and thermal content of biomass feedstocks. We use this case study to validate the models proposed and conduct numerical experiments. 

\section{Literature Review}\vspace{-0.08in}
The research presented in this paper is related to the following three main streams of literature: biomass blending problem, chance constraint programming and bilevel optimization.  

{\bf Biomass Blending Problem}: The research on biomass blending is scarce. The existing literature is mainly confined to the study of the impact that biomass blending has on the conversion performance \citep{shi2013impact}. These studies use sensitivity analysis to capture  the impacts of biomass supply/quality on conversion rate \citep{jacobson2014feedstock}. There exists literature about mathematical models that focus on optimizing blending of coal and grain products. For example,  \cite{sivaraman2002general} propose a model that identifies blending ratios to maximizes revenues from sales of blended grain products. Work by \cite{shih1995} proposes a multi-objective optimization model to identify coal blends which minimize system wide costs and greenhouse gas (GHG) emissions. To the best of our knowledge, there are no papers which focus on the impacts of biomass blending on supply chain costs. 

{\bf Chance Constraint Programming}: Many applications in supply chain \citep{LR07}, production planning \citep{MP00}, energy systems  \citep{wang2012chance}, etc. use chance constraint programming (CCP) to model uncertainties.  These models ensure that the probability of meeting a requirement  is above a certain threshold level. These models are typically very difficult to solve \citep{BL97} for the following two reasons. First, the probability of meeting a certain constraint cannot be computed exactly due to the computational challenge of multidimensional integration. Second, the feasibility region defined by the chance constraints may not be convex \citep{kim2015guide}. The computational difficulties for solving CCPs motivated the development of approximate solution approaches.  There are two main approaches to solve CCP models. The first approach discretizes the corresponding probability distribution and  solves the corresponding combinatorial problem \citep{DPR00, LA08}. The second approach develops convex approximations of the chance constraints \citep{NS06}. 


This research uses the sample average approximation (SAA) method. The approximation is obtained by replacing the
actual distribution in a chance constraint by an empirical distribution corresponding to a random sample. The resulting deterministic equivalent model ensures that the number of unexpected ``failures'' in these independent samples is below the model thresholds. Variations of the SAA for chance constrained problems have been investigated in \citep{APH08, LA08}. The theoretical properties of SAA have been studied in \cite{pagnoncelli2009sample}. These studies provide the conditions for which an upper and a lower bound to the original CCP problem can be obtained.

{\bf Bilevel Optimization}: A bilevel optimization model is a mathematical model with an optimization problem in the constraints. This model  is a generalization of the Stackelberg game. Many applications in resource planning, financial planning, land-use planning etc. are modeled and solved using  bilevel optimization models \citep{LSZ06}. Bilevel optimization models are difficult to solve since the corresponding feasible region is  not convex. A special case is the bilevel optimization model where the inner optimization model is linear. In this case, the inner optimization model can be replaced by the corresponding KKT conditions. The corresponding single-level problem is a nonlinear program. A number of approaches have been developed to solved these nonlinear programs \citep{Bard84, Bard98}. On addition to the KKT based approaches, other solution methods have been developed, such as, descent methods, penalty function methods and trust-region methods \citep{SMD17}.  

Within the framework of any bilevel optimization, a leader's decision is influenced by the reaction of his follower(s). In a setting with multiple followers,  leader's decision is influenced not only by the decision of each follower, but also by the relationships that exist among these followers. Work by \cite{LSZ06}  identifies nine different kinds of relationships amongst followers by establishing a general framework for bilevel multi-follower decision problems. One of the problems analyzed in great details is the \emph{uncooperative decision problem}, which is the same problem we solve in this research. \cite{LSZ06}  extend the Kuhn$-$Tucker approach to find an optimal solution for the uncooperative decision model and illustrate its performance via a real life case study.

\section{Blending Problem in a Centralized Supply Chain}
\subsection{Assumptions and Parameters}
Consider a supply chain with $|I|$  suppliers who provide biomass feedstock to a single biorefinery. The following parameters describe the biorefinery and are a function of its production capacity: $\alpha$ denotes the allowable ash content (in \%) and $\tau$ denotes the annual thermal requirement (in BTU). 

Let ${S}_{ib}$ denote the amount of biomass feedstock $b$ available at supplier $i$. Let  $\tilde{a}_{ib}$  denote the corresponding ash content, and $\tilde{h}_{ib}$ denote the thermal content for biomass type $b \in \mathcal{B}$. Ash and thermal contents are biomass qualities which impact the performance of the conversion process. We assume that ash  and thermal contents are stochastic and are represented by random variabels which follow some continuous distributions. 

\begin{figure}[htp]
	\centering
	\begin{subfigure}{0.45\textwidth}
		\includegraphics[scale=0.33]{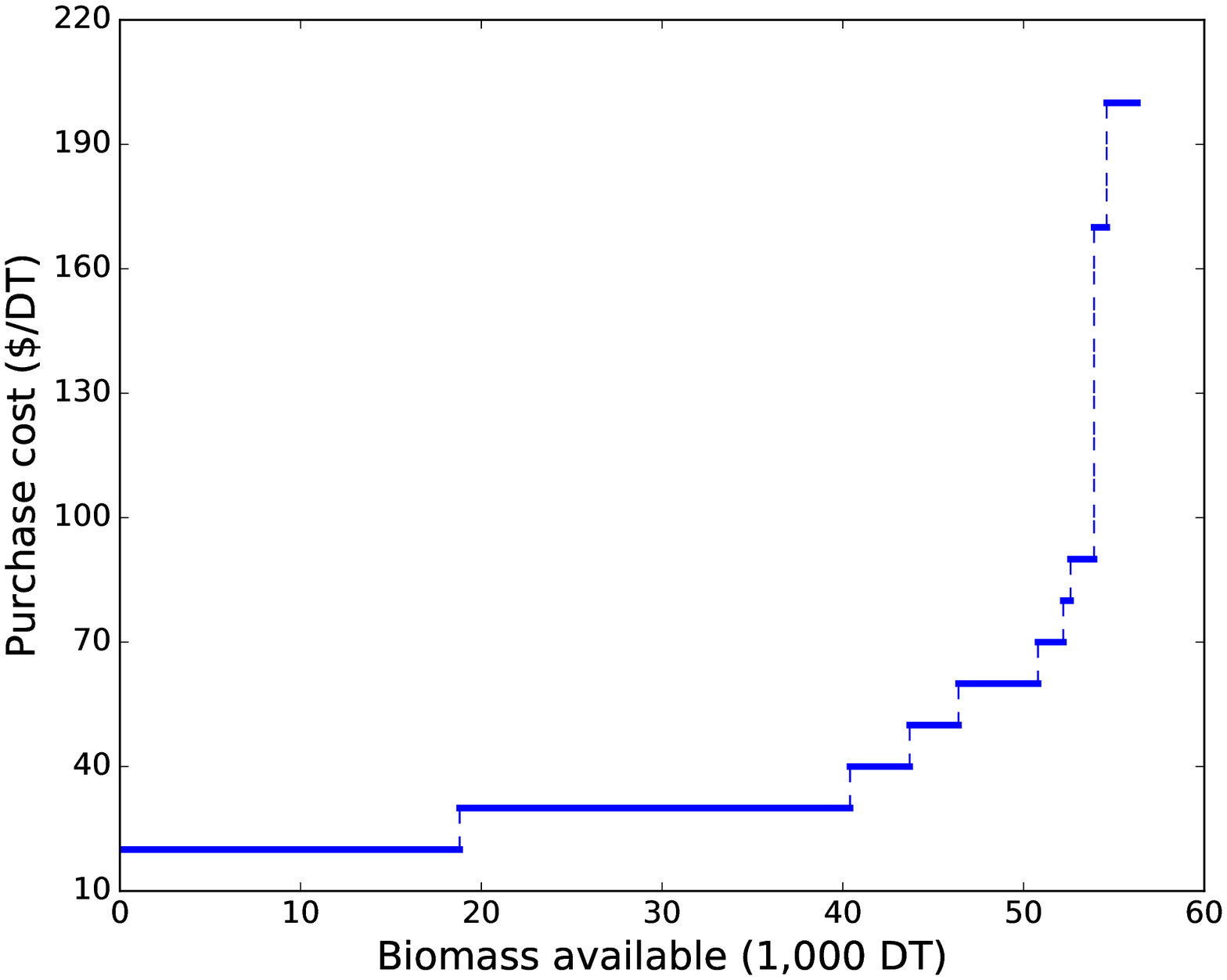}
		\caption{Supply Curve}
		\label{supplycurvea}
	\end{subfigure}
	\begin{subfigure}{0.45\textwidth}
		\includegraphics[scale=0.33]{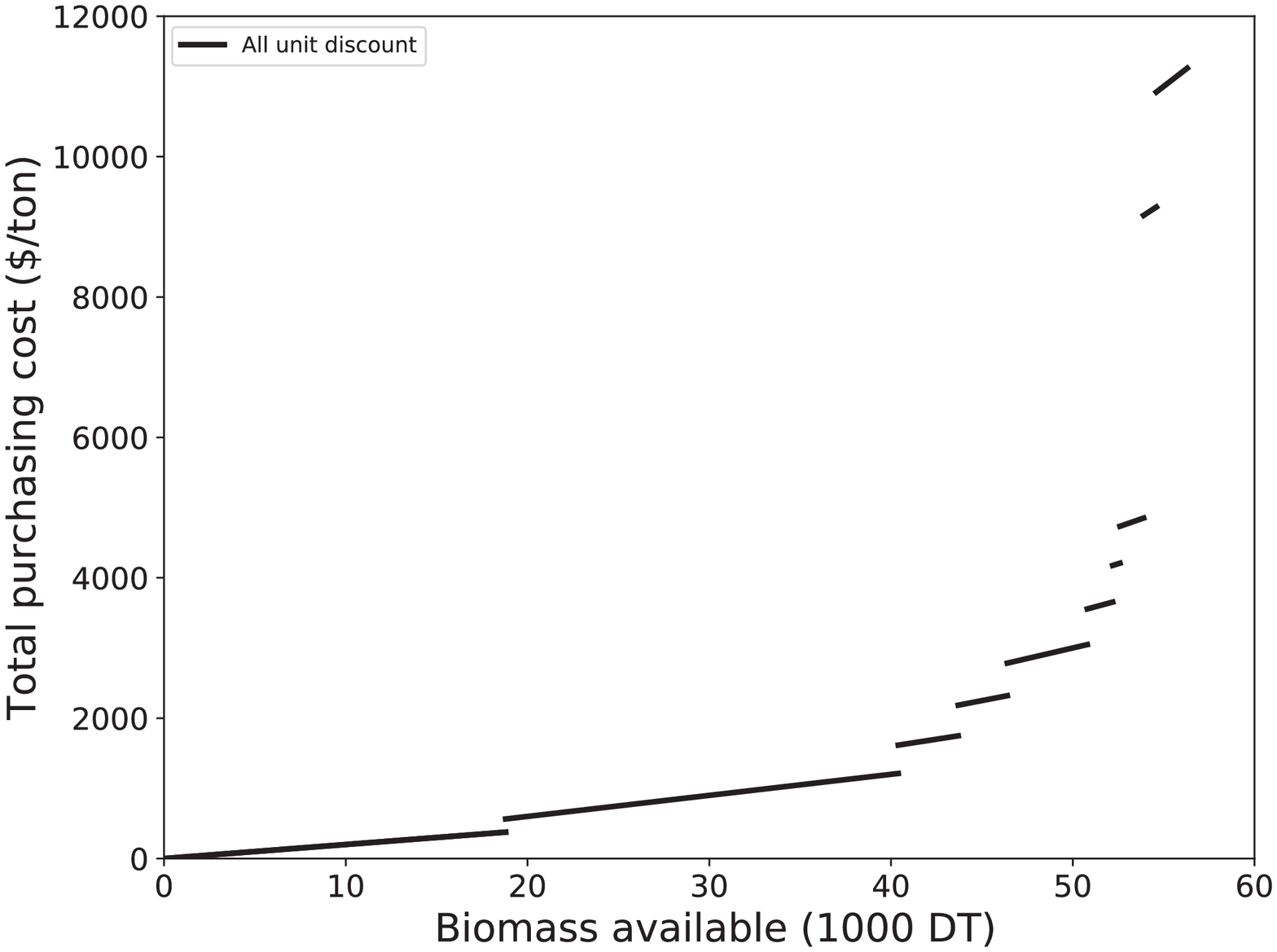}
		\caption{Total Purchasing Cost}
		\label{costcurvea}
	\end{subfigure}
	\caption{Supply Curve for Woody Biomass in Pickens County, SC.}
	\label{pricecurvea}
\end{figure}

The unit (farmgate) cost of biomass depends on the amount of biomass available. This cost increases with quantity. For example, wood residues is a type biomass and its harvesting and collection costs  impact the selling price. Two sources of woody residues are municipal waste (MW) from pruning trees in our backyards and  forest thinning. Harvesting cost of MW is zero, however, the amount of biomass available from MW is low. If a biorefinery is willing to pay the high price of  forest thinning, then additional amounts of biomass become available. The Billion Ton Study  conducted by the Oak Ridge National Laboratory (\citeyear{langholtz2016}) provides data about the amount of biomass available for a given set of {farmgate}  cost  for each county in the USA (see Figure \ref{supplycurvea}).
 Let $\mathcal{P}=\{10, 20, ...\}$ be the set of unit farmgate costs listed in the Billion Ton Study.  We use $p = 1,\ldots, |\mathcal{P}|$ to denote the index and $c_{p}$ to denote an element of this set. Notice that, the set of farmgate costs is the same for every biomass type. However, the amount of biomass available at this cost differs by supplier and biomass type (read Section 4.5 in \cite{langholtz2016}). 
Let $\underline{k}_{ibp}$ denote  the least and $\overline{k}_{ibp}$ denote the maximum amount of biomass available at cost $c_{p} \in \mathcal{P}$, thus,  ${S}_{ib} = \overline{k}_{ib|\mathcal{P}|}$.  
The relationship between  the amount available and the total purchasing cost at a biorefinery  is represented via a piece-wise linear function (see  Figure \ref{costcurvea}). 

Another problem parameter is $t_{ib} (= v_b\times Dist_i + g_b$) which denotes the unit transportation cost (in \$/ton). This cost depends on the type of biomass delivered and the distance from the supplier $i$, $Dist_i $; the variable unit cost, $v_b$; and  the fixed  unit transportation cost, $g_b$. Finally, $f_{b}$ denotes the cost of  processing and inventory  at the biorefinery. 

\subsection{A Stochastic Problem Formulation}
Let $X_{ib}$ denote the amount of biomass type $b$ purchased from supplier $i$. Then, the corresponding biomass purchase cost is represented by:


\begin{equation}\label{Function_AllUnitsDiscount}
\mathcal{F}_{ib}(X_{ib}) =\left\{ \begin{array}{lcl}
c_1 X_{ib}, & & \mbox{if } 0 \leq  X_{ib} \leq  \overline{k}_{ib1},\\
c_2 X_{ib},  & & \mbox{if } \underline{k}_{ib2} < X_{ib} \leq  \overline{k}_{ib2},\\
\\\ldots,\\
c_{|\mathcal{P}|} X_{ib},  & & \mbox{if } \underline{k}_{ib|\mathcal{P}|} < X_{ib} \leq  \overline{k}_{ib|\mathcal{P}|}. \\
 \end{array} \right.
\end{equation}
where, $\overline{k}_{ib,p} = \underline{k}_{ibp+1}$ for $p=1,\ldots, |\mathcal{P}|-1.$

The proposed  model  identifies a blendstock that minimizes the total supply chain costs including biomass purchasing, transportation, and processing and inventory costs at the biorefinery.
\begin{equation}
\begin{aligned}
    &\label{obj_s_1}\min:  \sum_{i \in I}\sum_{b \in B} \left[\mathcal{F}_{ib}(X_{ib}) + t_{ib}X_{ib} + f_bX_{ib}\right].
\end{aligned}
\end{equation}

This objective is minimized subject to the following constraints. Constraints (\ref{demand_1}) indicate that the amount of biomass shipped from supplier $i$ is limited by its availability.
\begin{equation}
\begin{aligned}
   \label{demand_1} X_{ib}  \leq S_{ib} \hspace{0.4in} \forall i \in I, b \in B.
 \end{aligned}
 \end{equation}

Since biomass quality impacts the performance of the conversion process, biorefineries require that the total ash content  be at most $\alpha$\% of biomass purchased, and the corresponding total thermal value be at least $\tau$ BTUs. However, these  are  soft requirements. That means, the biorefinery would like to meet these requirements. However, on a few occasions, the biorefinery is willing to violate these requirements if doing so will sufficiently decrease supply chain costs.

Let $\beta$ and $\gamma$ represent the risk parameters chosen by the biorefinery. These values are typically 10\% to 20\%. The following chance (probabilisitic)  constraint  indicates that ash content of biomass purchased by the biorefinery should be smaller than a threshold level $\alpha$ at least $(1-\beta)$ \% of the time.
\begin{equation} 
\begin{aligned}
  &\label{ash_L} Pr\left( \sum_{i \in I}\sum_{b \in B}\tilde{a}_{ib}X_{ib} \leq \alpha\sum_{i \in I}\sum_{b \in B}X_{ib}\right) \geq 1 - \beta.
\end{aligned}
\end{equation}

Similarly,  the following chance constraint indicates that the thermal content of biomass purchased by the biorefinery should be greater than the threshold level $\tau$ at least $(1-\gamma) \%$ of the time. Note that, the thermal energy gathered from the biomass delivered depends on its thermal content  ($\tilde{h}_{ib}$) and efficiency ($e_b$) of the conversion process.
\begin{equation} 
\begin{aligned}
    &\label{term_L} Pr\left(\sum_{i \in I}\sum_{b \in B} e_b \tilde{h}_{ib}X_{ib}\geq \tau\right) \geq 1- \gamma.
\end{aligned}
\end{equation}

(\ref{nonnegative_s}) are the non-negativity  constraints.
\begin{equation} 
\begin{aligned}
   \label{nonnegative_s}  X_{ib} \geq 0 \hspace{0.4in}\forall i \in I, b \in B.
\end{aligned}
\end{equation}

\noindent Formulation \eqref{obj_s_1} to \eqref{nonnegative_s}   is a chance constraint programming (CCP) model. We call this model  (${P}$).

\section{Solution Approaches: Centralized Blending Model}
Initially, we present a mixed-integer programming (MIP) formulation of model ($P$). Next, we propose a sample average approximation (SAA) of chance constraints \eqref{ash_L} and \eqref{term_L}. Finally, we present a linear approximation of $\mathcal{F}_{ib}(X_{ib})$ and a corresponding MIP formulation. 
\subsection{A Mixed-Integer Formulation of  (${P}$)}
Let $X_{ibp}$ be a decision variable which presents the amount of biomass type $b$ purchased from supplier $i$ which falls in bracket $p$. Let  $Z_{ibp}$ be a binary variable which takes the value 1 when the amount of biomass type $b$ purchased from supplier $i$ falls in bracket $p$, and takes the value 0 otherwise.

The following is a mixed-integer formulation of (${P}$).
\begin{subequations}\label{eq:MIP-1}
\begin{align}
(\bar{P}): \quad\quad & \label{obj_s_3}\min: \sum_{i \in I}\sum_{b \in B} \sum_{p\in \mathcal{P}} c_{ibp}X_{ibp}    \\
& \quad \text{s.t.} \nonumber \\
&\label{ash_s_3} Pr\left(\sum_{i \in I}\sum_{b \in B}\sum_{p\in \mathcal{P}}(a_{ib} - \alpha)X_{ibp} \leq 0\right) \geq 1 - \beta,\\
&\label{thermal_s_3} Pr\left(\sum_{i \in I}\sum_{b \in B} \sum_{p\in\mathcal{P}} e_bh_{ib}X_{ibp}\geq \tau\right) \geq 1- \gamma,\\
&\label{supply_cons} \quad \sum_{p\in \mathcal{P}} X_{ibp} \leq  S_{ib},\hspace{1.9in} \forall i \in I, b\in \mathcal{B},\\
&\label{resourcelimits_s_32} \underline{k}_{ibp}Z_{ibp}\leq X_{ibp} \leq  \overline{k}_{ibp}Z_{ibp}, \hspace{1.36in} \forall b \in B, i \in I, p \in \mathcal{P},\\
&\label{sumofZ} \quad \sum_{p\in \mathcal{P}}Z_{ibp} = 1, \hspace{2.05in} \forall i \in I, b\in \mathcal{B}, \\
& \label{nonnegative_s_3}  X_{ibp} \geq 0, \hspace{2.45in}\forall i \in I, b \in \mathcal{B}, p\in\mathcal{P}, \\
& \label{nonnegative_s_3a}  Z_{ibp} \in \{0,1\}, \hspace{2.21in}\forall i \in I, b \in \mathcal{B}, p\in\mathcal{P},
\end{align}
\end{subequations}

\noindent where, $c_{ibp} = c_p + t_{ib} + f_{b}.$ Notice that, constraints \eqref{supply_cons} are redundant because constraints  \eqref{resourcelimits_s_32} and \eqref{sumofZ}, and the definition of $\overline{k}_{ib|\mathcal{P}|} = S_{ib}$ impose the same restrictions on the solutions of (${P}$).

 \subsection{A Sample Approximation of Chance Constraints}\label{Modeling}
In order to make it easier for the reader to follow the description of the SAA, we provide the following succinct formulation of ($\bar{P}$). 
\begin{align} (\bar{P}): \hspace{0.6in} \min_{{x}\in \mathcal{Z}}  \,\, &Z = f({x}) \nonumber\\
 &p^1({x}, \tilde{a}) \leq {\beta}, \,\,\, \nonumber\\
 &p^2({x}, \tilde{h}) \leq {\gamma}.\nonumber
\end{align}

Let $E^1({x},\tilde{a}) = \sum_{i \in I}\sum_{b \in B}\sum_{p\in \mathcal{P}}(\tilde{a}_{ib} - \alpha){X}_{ibp}$, and, $p^1({x},\tilde{a}) = Pr \left( E^1({x},\tilde{a})>0 \right)$ (note: $Pr\left(E^1({x},\tilde{a})> 0\right) \leq \beta$, is equivalent to $Pr\left(E^1({x},\tilde{a})\leq 0\right) \geq 1-\beta$). Let $E^2({x},\tilde{h}) = \tau - \sum_{i \in I}\sum_{b \in B} \sum_{p\in \mathcal{P}}(e_b\tilde{h}_{ib}{X}_{ibp}$), and $p^2({x},\tilde{h}) = Pr\left(E^2({x},\tilde{h})  > 0\right)$. Let ${\vartheta}$ and $\mathcal{X}$ represent the optimal objective function value and the feasible region of ($\bar{P}$), respectively. Let $\mathcal{Z}$ represents solutions which satisfy constraints \eqref{supply_cons} to \eqref{nonnegative_s_3a}. We assume that, ($i$)  $\mathcal{X}$ is non-empty, and ($ii$) ${\vartheta}$  is bounded. 

The literature uses triangular and uniform  distributions to model $\tilde{a}_{ib}$ and $\tilde{h}_{ib}$, which are the random parameter in  ($\bar{P}$) \citep{SS16}. Thus, $E^1({x},\tilde{a})$ is the linear combination of $n_1$ ($n_1 = |I|\times |B| \times |\mathcal{P}|$ ) triangular distributed random variables; and $E^2({x},\tilde{h})$ is the linear combination of $n_2$ ($n_2 = |B| \times |\mathcal{P}|$) uniform distributed random variables.  We use Monte Carlo simulation to generate $N$ random samples of $\tilde{a}_{ib}$ and $\tilde{h}_{ib}$ from the corresponding distributions and use these values to calculate $E^1({x},\tilde{a})$, $E^2({x},\tilde{h})$. 


Let $P^1_N({x})=N^{-1}\sum_{s=1}^N\Delta(E^1_s({x},{a}_s))$ denote an empirical measure of the probability distribution of $E^1({x},\tilde{a})$  and $P^2_N({x})=N^{-1}\sum_{s=1}^N\Delta(E^2_s({x}, h_s))$ denote an empirical measure of the probability distribution of $E^2_s({x}, \tilde{h})$. Here, $a_s, h_s$, $s=1,\ldots, N$ are $N$ iid realizations of $\tilde{a}_{ib}, \tilde{h}_{ib}$, $\Delta()$ is a measure of probability mass function value, and $1/N$ is the probability assigned to each realization of  $\tilde{a}_{ib}, \tilde{h}_{ib}$. The SAA replaces the original distributions  of $E^1({x},\tilde{a})$, $E^2({x},\tilde{h})$ with $P_N^1({x})$,  $P^2_N({x})$ respectively \citep{PAS09}.  


Let $\mathbbm{1}_{(0,\infty)}: \mathcal{R}\rightarrow\{0,1\}$ be the indicator function of (0,$\infty$), i.e. 
\[\mathbbm{1}_{(0,\infty)}(t):=\left\{ \begin{array}{lrl}
1 & \text{if}&t> 0,\\
0 & \text{if}&t\leq 0.
\end{array}\right.
\]
We can now  approximate  $p^1({x}, \tilde{a})$ and $p^2({x}, \tilde{h})$ using the empirical measures $P^1_N$ and $P^2_N$ as follows
\begin{equation}\label{SAA_1p}\hat{p}_N^1({x}) = \mathbb{E}_{P_N^1}\left[\mathbbm{1}_{(0,\infty)}({E}^1({x},\tilde{a}))\right]=\frac{1}{N}\sum_{s=1}^N\mathbbm{1}_{(0,\infty)}({E}^1({x},a_s)),\end{equation}
\begin{equation}\label{SAA_2p}\hat{p}_N^2({x}) = \mathbb{E}_{P_N^2}\left[\mathbbm{1}_{(0,\infty)}({E}^2({x}, \tilde{h}))\right]=\frac{1}{N}\sum_{s=1}^N\mathbbm{1}_{(0,\infty)}({E}^2({x},h_s)).\end{equation}
Approximation $\hat{p}^1_N({x})$ returns what proportion of times  ${E}^1({x},\tilde{a})< 0$, and  $\hat{p}^2_N({x})$ returns what proportion of times  ${E}^2({x},\tilde{h})< 0$. The resulting SAA approximation model is presented below.   
\begin{eqnarray}  (\hat{P}): \hspace{0.6in} \min_{x \in\mathcal{Z}}  \,\, Z = f({x}) \nonumber \\
\label{SAA_1hat}  \hat{p}^1_N({x}) \leq \hat{\beta},\\
\label{SAA_2hat}  \hat{p}^2_N({x}) \leq \hat{\gamma}.
\end{eqnarray}
Let $\vartheta^N$ and $\mathcal{X}^N$ represent the optimal objective function value and the feasible region of ($\hat{P}$), respectively.  In this formulation, the reliability levels $1 -\hat{\beta}$ (for  $\hat{\beta}>0$) and $1-\hat{\gamma}$ (for $\hat{\gamma} > 0$) are different from the reliability level $1-\beta$ and $1-\gamma$ of the true model ($P$). 
Based on Theorem 5 in \cite{LA08}, if $\hat{\beta} < \beta$ and $\hat{\gamma} < \gamma$, every feasible solution of ($\hat{P}$) is feasible to (${P}$)  with high probability as $N$ gets large.  That is:  \[\vartheta^N \to {\vartheta} \mbox{ and } \mathcal{X}^N \to \mathcal{X}\ \mbox{ w.p.1 as } N \to \infty.\]


Constraints (\ref{SAA_1hat}) and (\ref{SAA_2hat})  of ($\hat{P}$) use indicator functions. Since commercial solvers cannot handle such functions, we reformulate these constraints by introducing the following continuous variables $\mathcal{V}, \mathcal{W}, \mathcal{U}, \mathcal{J}$ which quantify the violation of these constraints. The following are the equivalent linear constraints.
\begin{align}
{E}^1({x},a_s)+ \mathcal{V}_s - \mathcal{W}_s = 0, \,\,& \forall  s=1,\ldots,N,\label{saa_mip_1}\\
{E}^2({x},h_s)+ \mathcal{U}_s - \mathcal{J}_s = 0, \,\,& \forall  s=1,\ldots,N,\label{saa_mip_2}\\
\mathcal{V}_s, \mathcal{W}_s, \mathcal{U}_s, \mathcal{J}_s \geq 0, \,\,&\forall  s=1,\ldots,N. \label{saa_mip_3}
\end{align} 
This reformulation  minimizes the cost of violating the chance constraints \citep{charnes1955optimal, ABDELAZIZ20071811, ABDELAZIZ20121}.  Thus, variables $\mathcal{W}, \mathcal{U}$ appear also in the objective function, as follows:  
\begin{align}   
(\bar{\bar{P}}): \hspace{0.6in} \min_{x\in\mathcal{Z}}  \,\, Z= f({x}) + \lambda\sum_{s=1}^N  \mathcal{W}_s +\mu \sum_{s=1}^N  \mathcal{J}_s \nonumber \\
 \mbox{s.t. Constraints} \hspace{0.3in}\eqref{saa_mip_1}-\eqref{saa_mip_3}. \hspace{0.3in}\nonumber 
\end{align} 
This problem is easier to solve as compared to ($\hat{P}$). 
 Note however that, parameters $\lambda$ and $\mu$ are not known in advance. The size of these parameters is problem specific. When these penalties are too high, the minimization sets  $\mathcal{W}_s = 0$, and $\mathcal{J}_s = 0$ for all $s=1,\ldots, N$. Consequently, ${E}^1({x},a_s) \leq 0$ and ${E}^2({x},h_s) \leq 0$ for all  $s=1,\ldots, N.$  Thus, we develop an algorithm which identifies the value of $\lambda$ so that $\mathcal{W}_s = 0$ in at least  $\lceil(1-\hat{\beta})N\rceil$ of the scenarios generated; and identifies the value of $\mu$ so that $\mathcal{J}_s = 0$ in at least  $\lceil(1-\hat{\gamma})N\rceil$ of the scenarios generated. The {\bf SAA Algorithm} in Appendix A is an iterative procedure which uses a binary search to identify the values of  $\lambda$ and $\mu$.    


\subsection{A Linear Approximation of the Objective Function in (${P}$)} \label{sec:LPRelax_Central}
Let $F_{ib}(X_{ib})$ be a function defined as follows:
\begin{equation}\label{Function_AllUnitsDiscount_a}
{F}_{ib}(X_{ib}) =\left\{ \begin{array}{lcl}
{c}_{1}X_{ib}, & & \mbox{if } 0 \leq  X_{ib} \leq  \overline{k}_{ib1},\\
\lambda_{ib2} + {c}_{2}(X_{ib} - \underline{k}_{ib2}),  & & \mbox{if } \underline{k}_{ib2} < X_{ib} \leq  \overline{k}_{ib2},\\
\ldots \\
\lambda_{ib|\mathcal{P}|} + {c}_{|\mathcal{P}|}(X_{ib} - \underline{k}_{ib|\mathcal{P}|}), & & \mbox{if } \underline{k}_{ib|\mathcal{P}|} < X_{ib} \leq  \overline{k}_{ib|\mathcal{P}|}.\\
\end{array} \right.
\end{equation}
Where $\lambda_{ib1} = 0$ and $\lambda_{ibp} = \sum_{j \leq p-1} {c}_{j}(\overline{k}_{ibj} - \underline{k}_{ibj})$ for $p = 2,\ldots, \mathcal{P}.$ Function $F_{ib}(X_{ib})$ is a continues and convex approximation of function $\mathcal{F}_{ib}(X_{ib})$. Furthermore, $F_{ib}(X_{ib})$ provides an outer-approximation of $\mathcal{F}_{ib}(X_{ib})$  

The following is an approximations of ($P$).
\begin{eqnarray}
&({\mathbb{P}}): \hspace{0.6in} & \min: \sum_{i \in I}\sum_{b \in B} {F}_{ib} + \lambda \sum_{s=1}^N\mathcal{W}_s + \mu \sum_{s=1}^N\mathcal{J}_s\nonumber \\
&& \quad \text{s.t. } \nonumber\\
&& \quad  \eqref{saa_mip_1}-\eqref{saa_mip_3}, \\
&& \quad \label{pc_2} {F}_{ib} \geq \lambda_{ibp} + c_{p}(X_{ib}-\underline{k}_{ibp}), \hspace{0.8in} \,\,\,  \forall i\in I, b \in \mathcal{B}, p \in \mathcal{P},\\
&& \quad  \label{nonneg_cons}{X}_{ib} \geq 0, \hspace{2in} \,\,\,  \forall i\in I, b \in \mathcal{B}.
\end{eqnarray}

\begin{proposition}\label{prop1}
The feasible region of (${\mathbb{P}}$)  is convex. (Proof in Appendix B.) 
\end{proposition}
  
\begin{proposition}\label{prop2}
For each feasible solution of ($\bar{\bar{P}}$) one can find a feasible solution of (${\mathbb{P}}$), and vice versa. (Proof in Appendix B.) 
\end{proposition}
 
\begin{proposition}\label{Proposition2}
An optimal solution of ($\mathbb{P}$) is a lower bound of ($\bar{\bar{P}}$). (Proof in Appendix B.) 
\end{proposition}

Based on Propositions \ref{prop2} and \ref{Proposition2}, we develop an \emph{Algorithm for the  Centralized Problem} ($\bar{\bar{P}}$). The algorithm  solves (${\mathbb{P}}$) to obtain a feasible solution ${X}^*$.  Based on Proposition \ref{prop2},  ${X}^*$ is feasible for ($\bar{\bar{P}}$) as well. Next, we find an upper bound for ($\bar{\bar{P}}$) by calculating its objective function value at  ${X}^*$. In order to evaluate the quality of the approximation  we find a lower bound for ($\bar{\bar{P}}$) by calculating the objective function value of (${\mathbb{P}}$) at ${X}^*$. We report the corresponding error gap.

Note that,  \emph{Algorithm for the Centralized Problem}  solves ($\bar{\bar{P}}$) for a given value of $\lambda$ and $\mu.$ We use this algorithm within the {\bf SAA Algorithm} in order to identify the best values of $\lambda$ and $\mu$ that optimize ($\bar{\bar{P}}$). 

\section{Blending Problem in a Decentralized Supply Chain}
Due to the computational challenges of solving models for decentralized supply chains, most of the works in the literature assume centralized problem setting. This assumption makes the problem easy to solve, but often such a setting is not realistic. In this section we model this blending problem in a decentralized supply chain.  

We considers that the biorefinery and suppliers are independent entities  who have their own goals and objectives.  We propose a Stackelberg game  to model these relationships. The biorefinery is the leader of the game since it is  typically a large enterprise. Suppliers are the followers in the game since  farms in the USA are typically of small and medium size. Transportation costs in this supply chain are  high because biomass is bulk product and has low energy density. Therefore, to keep transportation costs low, a biorefinery  purchases from farms located nearby. As a result, we assume that one single biorefinery (the leader of the game) is located in the area. 

The game begins with the biorefinery setting  a ``door" price based on the type of biomass supplied. The ``door" price represents the amount of money paid to the supplier per ton of biomass delivered to the door of the biorefinery. The goal of the biorefinery is to identify a blendstock which minimizes its total supply chain costs while meeting thermal and ash content requirements. Next, suppliers decide how much to offer. If the amount offered meets the needs of the biorefinery, the game ends. Otherwise, the biorefinery adjusts the prices offered to suppliers, and the game continues. In this game, suppliers are independent, and each supplier focuses on maximizing his own profits.

 {\bf Leader's Problem:} The  biorefinery leads the game by setting a door price. Let $\mathcal{C}_{b}$ denote this price which is charged based on biomass type. 
 The objective of the leader is to identify a blendstock that minimizes her total supply chain costs. This objective is presented by the following equation. 
 
 \begin{equation} \nonumber 
\begin{aligned}
& \min: Z^{L} = \sum_{i \in I}\sum_{b \in \mathcal{B}} \sum_{p \in \mathcal{P}} (\mathcal{C}_{b} + f_{b})X_{ibp}. \\
\end{aligned}
 \end{equation}
Since biomass quality impacts the performance of the conversion process, the biorefinery requires that ash  content   \eqref{ash_s_3}, thermal content \eqref{thermal_s_3}, and non-negativity  \eqref{nonnegative_s_3} requirements are met. Additionally, we assume that farms will participate in this game only when it leads to profits. This assumption is realistic since, in a free market economy, farmers would not choose to participate in the game if doing so leads to an economic loss. Let $\bar{c}_{bp}$ represent the cost of harvesting, collecting and storing biomass $b$ within the cost bracket $p\in \mathcal{P}$.  $t_{ib}$ is the unit transportation cost.  We consider that each farm faces the same cost brackets which are defined in the centralized model ($\underline{k}_{ibp}, \overline{k}_{ibp}$).  Thus, the following is a constraint of leader's problem.
  \begin{equation} 
\begin{aligned}\label{profit_L}
& \qquad  \sum_{b \in \mathcal{B}}\sum_{p \in \mathcal{P}}(\mathcal{C}_{b} - \bar{c}_{bp} - t_{ib})X_{ibp} \geq 0,  \hspace{0.2in} \forall i \in I. \\
\end{aligned}
\end{equation}

Finally, 
\begin{equation}\label{L_nonneg}  
 \qquad \qquad \mathcal{C}_{b} \geq 0, \hspace{0.42in} b\in \mathcal{B}.
 \end{equation}
Let $\mathcal{X}^L$ denote the feasible region defined by constraints \eqref{ash_s_3}, \eqref{thermal_s_3},  \eqref{nonnegative_s_3},  \eqref{profit_L},  \eqref{L_nonneg}; and let the corresponding model be the \emph{leader's problem} ($P^L$).

{\bf Followers' Problem:} The objective of each farm is to identify how much biomass to supply in order to maximize profits. The objective of farm $i \in I$ is the following.
\begin{equation}\nonumber 
\begin{aligned}
& \qquad   \label{obj_s_4}\max: Z^{F}_{i} = \sum_{b \in \mathcal{B}}\sum_{p \in \mathcal{P}}(\mathcal{C}_{b} - \bar{c}_{bp} - t_{ib})X_{ibp}.
 \end{aligned}
\end{equation}

The amount of biomass delivered by a farm is impacted by its availability and  the price offered by the biorefinery.  Constraints \eqref{resourcelimits_s_32} to \eqref{nonnegative_s_3a} represent these restrictions.
 
Let $\mathcal{X}^F_i$ denote the feasible region of the problem faced by follower $i$  defined by  \eqref{resourcelimits_s_32} to \eqref{nonnegative_s_3a}. Let the corresponding  model be the \emph{followers' problem} ($Q^F_i$) and $\mathcal{Z}^F$ represent the set of optimal solutions. 
A complete formulation of the proposed  bilevel optimization model is presented in Appendix C. We refer to this as formulation ($Q$).  Next, we provide the corresponding succinct formulation in order to make it easier for the reader to follow the approach we develop for solving the bilevel optimization model proposed. 
\begin{align}  (Q): \hspace{0.6in} \min_{x^*, c^l}  \,\, &Z^L = F(x^*, c^l)  \nonumber \\
 & (x^*,c^l) \in \mathcal{X}^L \,\,\,  \nonumber\\
 & x^{*} \in \argmax_{x^f,z^f} Z^F_i = f_i(x^f_i, z^f_i) \hspace{0.3in}\forall i \in I  \nonumber\\
 &\hspace{0.4in}(x^f_i,z^f_i) \in \mathcal{X}^F_i. \,\,\, \nonumber
\end{align}

Let  ($\bar{Q}$) represent the SAA of model ($Q$) where the probabilistic constraints \eqref{ash_s_3} and \eqref{thermal_s_3}  are linearized as described in Section \ref{Modeling}.  A complete formulation of ($\bar{Q}$) is presented in Appendix C.

 \section{Solution Approaches: Decentralized Blending Model}\label{Sec:GAME}
We consider the optimistic (or strong) formulation of ($Q$). As such, the farms select from their set of optimal solutions, the amount of biomass to supply according to what is best for the biorefinery. 
  
In this section we focus on solving ($\bar{Q}$). Note that, ($\bar{Q}$) is not convex due to the inner optimization model (the followers' problem).  Also, ($\bar{Q}$) is nonlinear due to the bilinear terms $\mathcal{C}_{b}X_{ibp}$ in the objective. Next we present an exact solution, a heuristic  and a lower bound approach.  

 \subsection{An Exact Solution Approach to Solve ($\bar{Q}$)} 
 \subsubsection{Analyzing Followers' Problem:}\label{Sec:FPs}
Given the prices set by the leader  ($\mathcal{C}_{b}$), the followers' problem becomes an integer linear program. This problem can be decomposed by supplier and biomass type into $|I|\times|\mathcal{B}|$ sub-problems   of the following form. 
\begin{subequations}
\begin{align}
(\bar{Q}^F_{ib}): \hspace{0.6in}   \label{obj_s_7}& \max: Z^{F} =  \sum_{p \in \mathcal{P}} \overline{c}X_{p},\nonumber \\
 & \qquad \qquad \text{s.t.} \nonumber \\
& \qquad \qquad \sum_{p\in \mathcal{P}} X_{p} \leq  S,\\
& \qquad \qquad \label{cons_bounds}\underline{k}_{p}Z_{p}\leq X_{p} \leq  \overline{k}_{p}Z_{p}, \hspace{1in} \forall p \in \mathcal{P},\\
& \qquad \qquad \label{cons_Z1} \sum_{p \in \mathcal{P}} Z_{p} = 1,\\
& \qquad \qquad X_{p} \geq 0, \hspace{1.75in} \forall p \in \mathcal{P}, \\
& \qquad \qquad  \label{nonnegative_F1b} Z_{p} \in \{0,1\} \hspace{1.55in}\forall p \in \mathcal{P}.
\end{align}
\end{subequations}

\noindent Where, $\overline{c} = (\mathcal{C}_{b} - \overline{c}_{bp} - t_{ib})$. 
 \begin{proposition}\label{prop:lp}
The linear relaxation of $(\bar{Q}^F_{ib})$ provides an exact solution.  (Proof in Appendix B.)
 \end{proposition}
  \begin{lemma}
In an optimal solution to problem $(\bar{Q}^F_{ib})$, $\bar{c}\geq 0$ (due to constraints \eqref{profit_L}). Thus, at most one $X_{p}>0$ for $p\in\mathcal{P}$. Let $p^*$ be the cost bracket for which $X_{p}>0$, then,  $X_{p^*}=\overline{k}_{p^*}$ and  $X_p = 0$ for $p \in \mathcal{P} \setminus p^*$. (Proof in Appendix B.)
 \end{lemma}

\begin{theorem}
There is an $O(I\mathcal{B}\mathcal{P})$ algorithm which finds an optimal solution to the followers' problem ($\bar{Q}^F_{ib}$). 
\end{theorem}  

 {\bf Proof:}  For each supplier and biomass type, one can identify the  cost bracket which results in the maximum profits by following this two-steps procedure: (1) find $p^*=\argmax_{p\in\mathcal{P}} \bar{c}\overline{k}_p,$  (2) if $\bar{c}\overline{k}_{p^*}  < 0$, then,  $Z_p=X_p = 0$ for all $p\in\mathcal{P}$, (3) if $\bar{c}\overline{k}_{p^*} \geq 0$,  $X_{p^*}=\overline{k}_{p^*}$,  $Z_{p^*}=1$, and  $Z_p=X_p = 0$ for $p \in \mathcal{P} \setminus p^*$. We set $X_{p^*}=\overline{k}_{p^*}$, $Z_{p^*}=1$ when $\bar{c}\overline{k}_{p^*} = 0$ since we consider the optimistic formulation of ($P$). We call this procedure  {\bf Followers Algorithm}. 
 
 This two-steps procedure finds an optimal solution to $(\bar{Q}^F_{ib})$ in $O(\mathcal{P})$, which is the time it takes to find $p^*$ which leads to maximum profits. Since this procedure is followed by each supplier and for each biomass type, the total running time is $O(I\mathcal{B}\mathcal{Q})$. The corresponding solution found is optimal by construction. $\qed$

\subsubsection{A Single Level Optimization Model:}\label{Sec:KKT} 
Based on Proposition \ref{prop:lp}, the linear relaxation of  followers' problem provides an optimal solution. Thus, we replace followers' problem with its linear relaxation. The corresponding formulation is presented in Appendix C.  We reformulate the followers' problem by the corresponding KKT conditions to transform the bilevel optimization problem into a single level optimization model. The KKT conditions, which include the stationary constraints, primal feasibility constraints, dual feasibility constraints, and the complementary slackness constraints, are provided in Appendix C.  

The corresponding single-level model is not linear due to the bilinear terms $\mathcal{C}_{b}X_{ibp}$ in the objective and similar terms in the constraints. Thus, we solved the single-level model  using nonlinear solvers, such as, Couenne \citep{BeLeLiMaWa08} and POD \citep{nagarajan2016tightening, nagarajan2017adaptive}. POD uses an adaptive, multivariate partitioning of bilinear terms. POD is an iterative algorithm which exploits the advantages of piecewise polyhedral relaxation approaches via disjunctive formulations to solve nonlinear programs to global optimality. 

\subsection{A Heuristic Solution Approach to Solve ($\bar{Q}$)}\label{Sec:Heuristic} 
We propose an iterative procedure to solve ($\bar{Q}$).  This procedure starts by initializing the prices set by the leader  to $\hat{\mathcal{C}}_{b} = \min_{i \in \bar{\mathcal{I}}_b, p =1}\{t_{ib} + \overline{c}_{bp}\}$   for each $b \in \mathcal{B}$.  Here, $\bar{\mathcal{I}}_b := I$. Given these prices, we solve followers problems ($\bar{Q}_{ib}$) to identify how much should each follower supply to maximize his profits. Let $\bar{x}^f$ represent these quantities.  At the initial step, only farm $i^*_b (i^*_b= \argmin_{i \in \bar{\mathcal{I}}, p=1}\{t_{ib} + \overline{c}_{bp}\}$) would be offering biomass $b$  to the biorefinery. Next, the leader solves her  problem ($\tilde{Q}^L$) to identify whether the  quantities provided by the followers satisfy her needs. If the leader's problem is infeasible, then, she increases the prices offered. The prices are  increased even when the leader's problem is feasible in an effort to  attract suppliers that provide products of higher quality.   To find the new price, let  $\bar{\mathcal{I}}_b = \bar{\mathcal{I}}_b \setminus i^*_b$ and calculate  $\hat{\mathcal{C}}_{b} = \min_{i \in \bar{\mathcal{I}}_b, p =1}\{t_{ib} + \overline{c}_{bp}\}$   for each $b \in \mathcal{B}$. If $\bar{\mathcal{I}}=\emptyset$, then, let $p=p+1$, $\bar{\mathcal{I}}_b := I$ and calculate $\hat{\mathcal{C}}_{b}$.  This procedure continues until  no better solution is found for a fixed number of iterations  $\nu$.

 \begin{subequations} \label{Leader:heuristic}
 \begin{align}  
({\tilde{Q}}^L): \hspace{0.6in} \min  \,\, \tilde{\tilde{Z}}^L=F(x, c^l) + \lambda\sum_{s=1}^N  \mathcal{W}_s +\mu \sum_{s=1}^N  \mathcal{J}_s \nonumber \\
 \mbox{s.t. Constraints} \hspace{0.3in}\eqref{saa_mip_1}-\eqref{saa_mip_3} \hspace{0.3in}\nonumber\\
  c^l \geq 0, x\in \Pi_{i \in I} \mathcal{X}^F_i.\hspace{0.3in}\nonumber\\
 x \leq \bar{x}^f. \hspace{0.7in}\,\,\, \nonumber
\end{align} 
%
%
\end{subequations}

  
\begin{proposition}
Solutions found by the proposed heuristic are feasible for model ($\bar{Q}$).
\end{proposition}
It is straightforward to see that the proposed heuristic generates feasible solutions since the heuristic stops when, at a given price offered by the leader, her problem is feasible; and the corresponding amounts of biomass offered  maximize farmers' profits.  

\subsection{A Lower Bound for ($\bar{Q}$)}\label{sec:Dec_Relax} 
Consider the following relaxation of model ($\bar{Q}$). 
\begin{align}  
(\tilde{\tilde{Q}}^L): \hspace{0.6in} \min  \,\, \tilde{\tilde{Z}}^L=F(x, c^l) + \lambda\sum_{s=1}^N  \mathcal{W}_s +\mu \sum_{s=1}^N  \mathcal{J}_s \nonumber \\
 \mbox{s.t. Constraints} \hspace{0.3in}\eqref{saa_mip_1}-\eqref{saa_mip_3} \hspace{0.3in}\nonumber\\
  c^l \geq 0, x\in \Pi_{i \in I} \mathcal{X}^F_i.\hspace{0.3in}\nonumber
\end{align} 

This model minimizes the objective of the leader. The corresponding feasible region is the intersection of the feasible regions of the followers' problems and the feasible region of the leader. Since this model formulation is a relaxation of ($\bar{Q}$), its optimal objective function value is a lower bound for  ($\bar{Q}$). We use the corresponding lower bound in order to evaluate the performance of the heuristic solution approach proposed in Section \ref{Sec:Heuristic}.

\section{Numerical Analysis}
\subsection{Case Study Description and Related Data}
The main source of data is the Billion Ton Study \citep{langholtz20162016}. We focus our study in South Carolina, and used the county-level data of biomass supply.  Table \ref{tab:InputData} lists the different types of biomass available in South Carolina which are suitable to  use in a thermochemical conversion process. For each type of biomass and different scenarios, the report presents the expected amount available during 2014-2040. We used the data corresponding to these scenarios: ``Medium housing, low energy demands" for forest biomass;  ``Wastes and other residues" for C\&D and MSW waste;  ``Base case, single energy crops" for hybrid poplar. For poplar, we use the data of 2026, which is the most recent data available in the report. For the rest of biomass types, we use the data of 2016.

Tables \ref{tab:InputData} and \ref{tab:refinery} summarize the input data we use in the numerical analysis. This data was collected from \cite{jacobson2014feedstock}, Bioenergy Feedstock Library \citep{BFL} and \cite{harris2004potential}. Table \ref{tab:InputData} presents the average ash content (AAC), the range of ash content before (ACR)  and after ($\overline{\mbox{ACR}}$) pre-processing, low heating value (LHV), and high heating value (HHV)  for different types of biomass; and, harvesting and collection (H\&C) cost, processing cost (Pr) storage cost (St), fixed ($g_b$) and variable ($v_b$) transportation costs. 

\begin{table}[htp!]
\caption{Summary of Input Data}
\scriptsize
\centering
{%
\setlength{\tabcolsep}{0.28em}
\begin{tabular}{lccccccccccc}
\toprule
  &\multicolumn{5}{c}{Physical Properties} & \multicolumn{5}{c}{Costs}\\
 \cmidrule(lr){2-6}
 \cmidrule(lr){7-11}
 &AAC&ACR&$\overline{\mbox{ ACR}}$&LHV&HHV&H\&C&Pr&St&$g_b$&$v_b$&TQ \\
Feedstock & ($ wt.\% $)& ($ wt.\% $)&($ wt. \% $)&($\frac{10{^6}BTU}{DT}$) & ($\frac{10{^6}BTU}{DT}$)
 & (\$/DT) & (\$/DT)& ( \$/DT) & (\$/DT) & (\$/DT/mile)& (MDT) \\
 \midrule
Hybrid Poplar &0.50&0.3 - 4.3&0.30 - 0.75&16.768 &16.982 &22.24&23.97&3.23&20.53&0.046&0.34\\
 Pine&0.75&0.1 - 6.0&0.10 - 1.13&14.510&15.656&20.19&12.85&3.23&20.53&0.046 & 0.60\\
 SP$^{*}$ Residue&1.00&0.8 - 2.2&0.80 - 1.50&15.232&17.202&0.00&23.97&3.23&20.69&0.046&0.11 \\
 SN$^{**}$ Residue &1.00&0.8 - 2.2&0.80 - 1.50&15.232&17.202&0.00&23.97	&3.23&20.69&0.046& 0.27  \\
 Mixed Residue &1.20&0.8 - 2.2&0.80 - 1.80&15.160&17.892&0.00&23.97	&3.23&20.69&0.046& 0.25 \\
 C\&D$^{***}$ Waste &1.00&0.8 - 2.2&0.80 - 1.50&14.510&17.648&0.00&28.12&3.23&22.87&0.046& 0.34 \\
 MSW$^{****}$ Waste &10.00&7.0 - 15.0&7.00 - 15.00&10.250&13.680&0.00&19.70	&4.50&20.69&0.046& 0.099 \\
 \bottomrule
\multicolumn{12}{l}{$^{*}$ softwood planted, $^{**}$ softwood natural, $^{***}$ construction and demolition, $^{****}$ municipal solid waste}
\end{tabular}}
\label{tab:InputData}
\end{table}%
In our numerical analysis we conduct sensitivity analysis with respect to biorefinery capacity. Table \ref{tab:refinery} presents biorefinery capacities in terms of thermal requirement and biomass supply.  

\begin{table}[htp!]
    \centering
    \scriptsize
    \caption{Biomass Refinery Requirements}
    \begin{tabular}{ccccc}
    \toprule
     Demand & Allowable Ash Content & Thermal Requirement &Thermal Efficiency  \\
     (MDT/year) &($\%$) & ($10^9$ BTU/year) &($\%$) \\
    \midrule
    $0.3-0.6$ & $\leq1$ & $3,838-7,677$& $75-80$ \\
    \bottomrule
     \end{tabular}
    \label{tab:refinery}
\end{table}%

We consider each county to be a biomass supplier, and we assume that the available biomass is located at the centroid of this county. We consider that ash content for each supplier follows a triangular distribution   with  mean and range as defined in Table \ref{tab:InputData}. We consider that biomass heating value, for each supplier and biomass type, is uniformly distributed with bounds defined by LHV and HHV presented in Table \ref{tab:InputData} \citep{SS16}. Our model does not consider facility location decisions. Thus, we identify a biorefinery location in a separate model which minimizes the weighted travel distance to all suppliers.

The algorithms proposed are  programmed in Julia 0.6.2 using modeling language JuMP. These models ran on Clemson University's high performance resource Palmetto Cluster and used 8 nodes and 64 GB RAM. The linear and mixed integer  programs are solved using GUROBI callable subroutines. 
\subsection{Evaluating the Performance of the Solution Approaches Proposed}
\subsubsection{Linear Approximation Model ($\mathbb{P}$):}

Table \ref{tab:MIPvsLPRelax} presents the average, minimum and maximum gap between the feasible solutions  and the corresponding lower bounds of ($\bar{\bar{P}}$) found by solving ($\mathbb{P}$) via the \emph{Algorithm for the Centralized Problem}  proposed in Section \ref{sec:LPRelax_Central}.  
The average (over 10 replications) error gap is less than 0.14\% and the running time is smaller than 0.25 seconds which demonstrate that the proposed approximation can provide high quality solutions in a short amount of time.   

\begin{table}[htp!]
	\caption{Summary of the Results of Linear Approximation Model}
	\label{tab:MIPvsLPRelax}
	\scriptsize
	\centering	
	\begin{tabular}{cccc}
		\toprule
		& \multicolumn{3}{c}{Error Gap (in  $\%$)}                                                                                                                                                                                                                   \\
		\cline{2-4}
		 Demand&Avg&Min&Max\\
		\cmidrule(lr){1-1}
		\cmidrule(lr){2-2}
		\cmidrule(lr){3-3}
		\cmidrule(lr){4-4}
		0.3&0.02&0.00&0.07\\
		0.4&0.02&0.00&0.07\\
		0.5&0.03&0.00&0.08\\
		0.6&0.07&0.02&0.18\\
		0.7&0.09&0.02&0.24\\
		0.8&0.14&0.04&0.26\\
		\bottomrule	
	\end{tabular}
\end{table}

\subsubsection{Single-Level Optimization Model to Solve ($\bar{Q}$):}  
The single-level optimization model is not linear due to a few nonlinear terms in the objective and constraints. In order to solve this  Mixed-Integer Nonlinear Program (MINLP) we used  Couenne and POD. Both solvers failed to solve instances with more than 2 suppliers, 1 biomass type and 2 cost backets. Thus, we created a small size problem instance with 2 suppliers who  supply  hybrid poplar. We picked hybrid poplar since it has low ash content. For this problem we only generated 2 scenarios. The results are summarized in Table \ref{tab:SmallSizeCompare}.  

\begin{table}[htp!]
	\caption{Evaluating the Single-Level Optimization Model for A Small Size Problem ($\beta = 0.1$, $\gamma = 0.1$)}
	\scriptsize
	\centering
	{%
		\begin{tabular}{cccccc}
			\toprule
			& &&\multicolumn{2}{c}{Violation (\%)}&\\
			\cmidrule(lr){4-5}
			MINLP &Demand&Costs&& &Running\\
			Solver &(MDT/year)&(\$/DT)  &Ash &Thermal &Time (sec.)\\
			\midrule
			Couenne & 0.3&128.59&0.0&0.0&267.11\\
			\midrule
			POD &0.3&133.48&0.0&0.0&181.83\\
			\midrule
			Heuristic & 0.3&126.00&0.0&0.0&6.45\\
			\bottomrule
	\end{tabular}}
	\label{tab:SmallSizeCompare}
\end{table}%

The result indicate that the running time of the proposed Heuristic (see Section \ref{Sec:Heuristic})  is order of magnitude smaller than the running time of Couenne and POD. The  Heuristic also provides a feasible solution of  higher quality which is 2\% lower than the solution found from Couenne and 6\% lower than the solution found by POD. Notice that, we use the ``trunk" version of Couenne solver, thus the solution found is not guaranteed to be optimal (see \cite{Couenne}), which is the case with the problem solved.

\subsubsection{Heuristic Algorithm to Solve ($\bar{Q}$):}  
In order to evaluate the performance of the proposed Heuristic in solving large problem instances, we  compare its solutions with the lower bounds found from solving  ($\tilde{\tilde{Q}}^L$) (see Section \ref{sec:Dec_Relax}). 
 ($\tilde{\tilde{Q}}^L$) is a bilinear program which we solve using Couenne. Since the time it takes to solve this problem in  Couenne is too long, we only solved the following problems. Problem 1 considers the whole dataset. Problem 2 considers a smaller dataset consisting only of suppliers of hybrid poplar and softwood residues. Problem 3 considers a supply chain with 23 suppliers (rather than the 46 suppliers we have in our dataset). We present the objective function value obtained from solving these problems and the corresponding error gap. 
 
 When solving Problem 1, we stoped Couenne after 48 hours. The solution found is not optimal. For Problems 2 and 3, we stoped Couenne after 10 hours. The solutions found for both problems are not optimal. Thus, we cannot claim that the objective function values of  ($\tilde{\tilde{Q}}^L$) are valid lower bounds for the Heuristic. However, the  objective function values of Problems 1 and 2 are within 1\% of the objective function value of the Heuristic.
  
\begin{table}[htp!]
	\caption{Comparison of ($\tilde{\tilde{Q}}^L$) formulation and the Heuristic Algorithm (N = 1) }
	\scriptsize
	\centering
	{%
		\begin{tabular}{cccc|ccc|ccc}
			\toprule
			\multicolumn{4}{c|}{\bf Problem 1} & \multicolumn{3}{c|}{\bf Problem 2} &  \multicolumn{3}{c}{\bf Problem 3}\\
			\hline
			&($\tilde{\tilde{P}}^L$)&Heuristic & Error Gap & ($\tilde{\tilde{P}}^L$)&Heuristic & Error Gap & ($\tilde{\tilde{P}}^L$)&Heuristic & Error Gap \\
			&& & (in \%) & & & (in \%)& & & (in \%)\\
			\cline{1-10}
			Obj. Func. & 29,633,223 &	 29,069,403  &  -1.900 & 29,999,309 &  30,003,718 & 0.015 & 40,829,002 &	 30,684,451 & -24.846\\
			Value  & \multicolumn{3}{c|}{} &\multicolumn{3}{c|}{}  & \multicolumn{3}{c}{}  \\
			Run Time& 172,834 &	61 & & 35,736 &	193 & & 35,754 &	0.013\\	
			(sec)& \multicolumn{3}{c|}{} &\multicolumn{3}{c|}{}  & \multicolumn{3}{c}{}  \\		
			\bottomrule
	\end{tabular}}
	\label{tab:Rel_Alg_Full}
\end{table}%

\subsection{Managerial Insights}
Tables \ref{tab:centrl_both_03} and \ref{tab:centrl_both_02} summarizes the results of solving ($P$) for different levels of demand.   
Based on these results, SN residues contributes  30 to 45\% of the amount in a blend and pine contributes  16 to 18\%, and SP residues contribute 13 to 25\% of the amount blended. MSW waste has not been utilized and the use of hybrid poplar increases with demand for biomass. Since hybrid poplar is expensive, this increase in utilization impacts the unit cost of the blend. 

Comparing the results of Tables \ref{tab:centrl_both_03} and \ref{tab:centrl_both_02}, one can observe that, as the risk level decrease (i.e. $\beta$ and $\gamma$ decrease), the chance constraints become more restrictive, thus, the  blends identified contain greater amounts pine and hybrid poplar since these biomass types have lowest ash contents.   The running time of the heuristic algorithm is higher when $\beta = \gamma = 0.2$ (compared to $\beta = \gamma = 0.3$) since the sample size $N$ is larger.


	\begin{table}[htp!]
		\caption{Costs and Biomass Blending Ratios Under Ash and Thermal Content Uncertainties  ($\beta =  \gamma$ = 0.3)}
		\scriptsize
		\centering
		{%
			\begin{tabular}{ccccccccccc}
				\toprule
				& &&\multicolumn{7}{c}{Blending Ratios for Thermochemical Conversion Process (in \%)}&\\
				\cmidrule(lr){4-10}
				Demand & Cost&Costs& Hybrid &&SP & SN & Mixed& C\&D &  &Running\\
				(MDT/year) &(\$1000)&(\$/DT) & Poplar&Pine&Residue&Residue& Residue& Waste& MSW &Time (sec.)\\
				\midrule

				0.3& 29,396	&87.76	&6.59	&18.54	&24.16	&45.11	&3.55	&2.05	&0.00	&0.20\\
				0.4& 39,462	&88.49	&8.11	&17.55	&21.31	&43.53	&5.53	&3.97	&0.00	&0.19\\
				0.5& 49,678	&89.31	&9.85	&16.29	&18.66	&39.98	&8.34	&6.87	&0.00	&0.19\\
				0.6& 60,039	&90.07	&10.97	&15.77	&16.34	&36.78	&10.86	&9.29	&0.00	&0.18\\
				0.7& 70,519	&90.70	&11.58	&15.62	&14.27	&33.33	&12.69	&12.52	&0.00	&0.18\\
				0.8& 81,123	&91.23	&11.47	&16.34	&12.54	&29.82	&14.27	&15.56	&0.00	&0.18\\
				\bottomrule
				\multicolumn{7}{l}{Note: $\hat{\beta} =\hat{\gamma}$ = 0.3 and $N = 125.$}
		\end{tabular}}
		\label{tab:centrl_both_03}
	\end{table}%
	
	\begin{table}[htp!]
		\caption{Costs and Biomass Blending Ratios Under Ash and Thermal Content Uncertainties ($\beta = \gamma$ = 0.2)}
		\scriptsize
		\centering
		{%
			\begin{tabular}{ccccccccccc}
				\toprule
				& &&\multicolumn{7}{c}{Blending Ratios for Thermochemical Conversion Process (in \%)}&\\
				\cmidrule(lr){4-10}
				Demand & Cost&Costs& Hybrid &&SP & SN & Mixed& C\&D &  &Running\\
				(MDT/year) &(\$1000)&(\$/DT) & Poplar&Pine&Residue&Residue& Residue& Waste& MSW &Time (sec.)\\
				\midrule
				0.3& 29,519& 	88.03&	6.73&	19.47&	24.42&	44.36&	2.47&	2.54&	0.00&	0.36\\
				0.4& 39,614& 	88.88&	9.19&	17.08&	21.76&	42.67&	5.11&	4.19&	0.00&	0.42\\
				0.5& 49,842& 	89.58&	10.53&	16.23&	18.86&	39.91&	7.47&	7.00&	0.00&	0.40\\
				0.6& 60,220& 	90.31&	11.65&	15.52&	16.51&	36.45&	9.81&	10.05&	0.00&	0.40\\
				0.7& 70,755& 	91.02&	12.49&	15.25&	14.24&	33.02&	12.16&	12.84&	0.00&	0.39\\
				0.8& 81,404& 	91.64&	13.01&	15.33&	12.53&	29.74&	14.04&	15.34&	0.00&	0.38\\
				\bottomrule
				\multicolumn{7}{l}{Note: $\hat{\beta} =\hat{\gamma}$ = 0.2, and  $N = 250.$}
		\end{tabular}}
		\label{tab:centrl_both_02}
	\end{table}%

 Tables \ref{tab:decentrl_both03} and \ref{tab:decentrl_both02} summarizes the results from solving ($\bar{Q}$) for different levels of demand.     
Via these experiments we find that SN residues  and pine count for about 50-60\% of the blends identified. This is mainly because pine has low ash content and SN residues are not expensive.  As demand increases, we observe an increase in the amount of  pine used. This increase impacts the cost of the blend. The amount of mixed residues and C\&D waste also increases with demand. This is mainly due to the low cost of delivering these biomass types, and the limited amount of SN residues available. MSW was not  used in a blend due to its high ash content. Hybrid poplar is used in moderation due to its high cost, although, its ash content is low. Based on these results, decreasing risk level (i.e., $\beta, \gamma$ decrease) leads to higher costs in the supply chain. 


\begin{table}[htp!]
		\caption{Costs and Biomass Blending Ratios Under Ash and Thermal Content Uncertainties ($\gamma$ = 0.3, $\beta$ = 0.3)}
		\scriptsize
		\centering
		{%
			\begin{tabular}{ccccccccccc}
				\toprule
				& &&\multicolumn{7}{c}{Blending Ratios for Thermochemical Conversion Process (in \%)}\\
				\cmidrule(lr){4-10}
				Demand & Cost&Cost& Hybrid &&SP & SN & Mixed& C\&D &  &Running\\
				(MDT/year) &(\$1000)&(\$/DT) & Poplar&Pine&Residue&Residue& Residue& Waste& MSW &Time (sec.)\\
				\midrule
				
				0.3& 30,463	&91.11	&7.61	&16.42	&26.05	&46.30	&2.05	&1.58	&0.00	&708\\
				0.4& 41,159	&93.39	&9.56	&19.93	&5.46	&44.49	&11.74	&10.19	&0.00	&773\\
				0.5& 51,651	&93.85	&7.99	&21.97	&6.40	&38.60	&12.21	&14.51	&0.00	&750\\
				0.6& 62,403	&94.46	&7.26	&22.98	&7.44	&34.05	&12.94	&17.08	&0.00	&777\\
				0.7& 73,305	&95.12	&7.55	&22.47	&8.41	&30.52	&13.21	&19.58	&0.00	&703\\
				0.8& 84,175	&95.58	&7.17	&23.32	&8.71	&27.56	&13.96	&21.14	&0.00	&696\\
				
				\bottomrule
				\multicolumn{7}{l}{Note, $\hat{\beta} = \hat{\gamma} = 0.0, N = 125.$}
		\end{tabular}}
		\label{tab:decentrl_both03}
	\end{table}%

	\begin{table}[htp!]
		\caption{Costs and Biomass Blending Ratios Under Ash and Thermal Content Uncertainties ($\gamma$ = 0.2, $\beta$ = 0.2)}
		\scriptsize
		\centering
		{%
			\begin{tabular}{ccccccccccc}
				\toprule
				& &&\multicolumn{7}{c}{Blending Ratios for Thermalchemical Conversion Process (in \%)}&\\
				\cmidrule(lr){4-10}
				Demand & Cost&Costs& Hybrid &&SP & SN & Mixed& C\&D &  &Running\\
				(MDT/year) &(\$1000)&(\$/DT) & Poplar&Pine&Residue&Residue& Residue& Waste& MSW &Time (sec.)\\
				\midrule

				0.3& 30,585	&91.54	&8.93	&15.82	&26.07	&46.18	&1.19	&1.80	&0.00	&1,940\\
				0.4& 41,096	&93.74	&9.78	&20.91	&6.27	&44.91	&10.27	&9.89	&0.00	&1,947\\
				0.5& 51,664	&94.33	&8.89	&22.16	&7.25	&39.35	&11.87	&12.68	&0.00	&1,920\\
				0.6& 62,543	&95.03	&9.05	&21.70	&8.42	&35.12	&11.96	&15.78	&0.00	&1,789\\
				0.7& 73,641	&95.67	&9.11	&21.83	&8.96	&31.45	&13.19	&17.22	&0.00	&1,729\\
				0.8& 84,841	&95.98	&7.17	&24.79	&9.43	&27.98	&13.25	&18.99	&0.00	&1,646\\
				
				\bottomrule
				\multicolumn{7}{l}{Note, $\hat{\beta} = \hat{\gamma} = 0.0, N = 250.$}
		\end{tabular}}
		\label{tab:decentrl_both02}
	\end{table}%

\begin{figure}[h]
	
	\centering
	
	\resizebox{1.02\textwidth}{!}{
		\begin{subfigure}{8.75cm}
			\begin{tikzpicture}[baseline]
			\captionsetup{type=figure,justification=centering}
			\begin{axis}[%
			xmin=0, xmax=7,%
			ymin=1, ymax=6,%
			xtick={1,2,3,4,5,6},xticklabels={0.3,0.4,0.5,0.6,0.7,0.8},%
			xlabel=Biorefinery Demand,
			ylabel= Gap ($\%$),
			]
			\boxplot{1}{3.6014464120128}{3.53943167857738}{3.65311073458992}{3.50648413126646}{3.78388340631539}
			\boxplot{2}{4.05039807235034}{2.75380798229604}{4.45467427959378}{1.79029024546853}{5.8168468567535}
			\boxplot{3}{3.61062024342534}{2.92269391174021}{4.11358959751433}{2.40271702989313}{5.14681724810486}
			\boxplot{4}{4.07714773640468}{3.2543901343435}{4.56611445533715}{1.83222975879606}{5.15619610621638}
			\boxplot{5}{4.15558047918597}{3.87760635325207}{4.38432929016348}{3.43394080442727}{4.4905577944185}
			\boxplot{6}{4.09601498658496}{3.861351685913}{4.52077944159822}{3.62142662813878}{5.21283327842164}
			\end{axis}
			\end{tikzpicture}
			\caption{$\beta$ = 0.2, $\gamma$ = 0.2}
			\label{fig:a}
		\end{subfigure}
		
		\begin{subfigure}{8.75cm}
			\begin{tikzpicture}[baseline]
			\captionsetup{type=figure,justification=centering}
			\begin{axis}[%
			xmin=0, xmax=7,%
			ymin=1, ymax=6,%
			xtick={1,2,3,4,5,6},xticklabels={0.3,0.4,0.5,0.6,0.7,0.8},%
			xlabel=Biorefinery Demand,
			ylabel= Gap ($\%$),
			]
			\boxplot{1}{3.65240193827649}{3.58233509050866}{3.68518286233624}{3.44530596921677}{3.75879417402833}
			\boxplot{2}{4.5268701357773}{4.21752456034392}{4.69911650353527}{1.45585352408782}{5.2652923340207}
			\boxplot{3}{4.01090002709982}{3.59778818827606}{4.36215428442101}{3.11635493347103}{4.97138637967642}
			\boxplot{4}{3.98677736254218}{3.45596451738323}{4.26114116031466}{2.74418265055396}{5.53850506454628}
			\boxplot{5}{3.76812886115485}{3.43624936385599}{4.47924858525377}{2.92645233975557}{5.30845972972324}
			\boxplot{6}{3.75975859244084}{3.28289405217453}{4.13693597107494}{2.56676610133422}{5.28359813306326}
			\end{axis}
			\end{tikzpicture}
			\caption{$\beta$ = 0.3, $\gamma$ = 0.3}
			\label{fig:b}
		\end{subfigure}

	}
	\caption{Decentralize vs Centralize - Percent  Gap}
	\label{fig:Dec_Cen_Gap}
\end{figure}
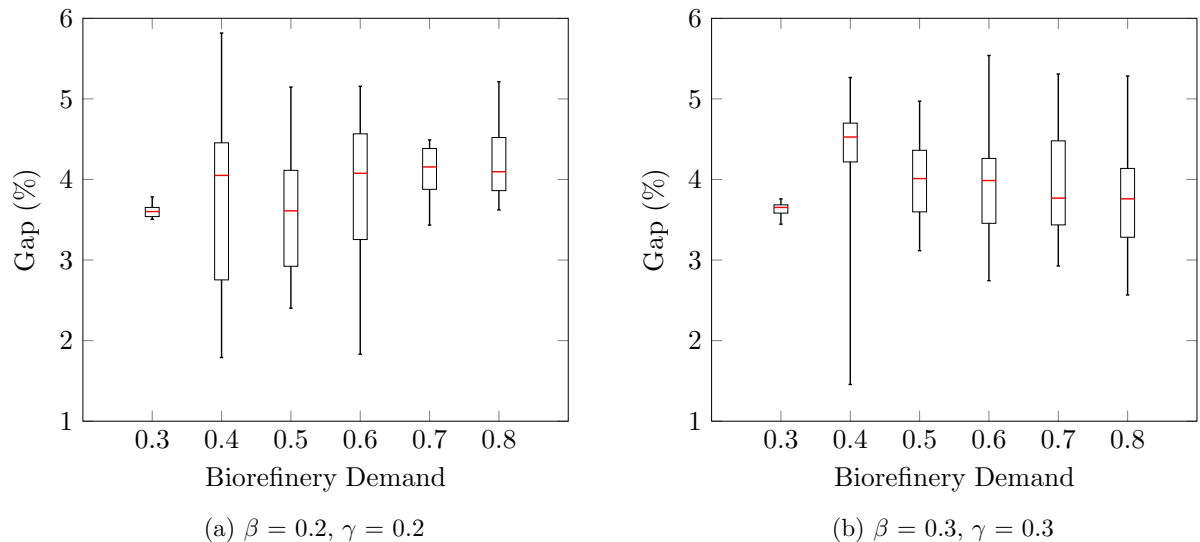

Figure \ref{fig:Dec_Cen_Gap} summarizes the gap between the objective function valued of the  centralized and decentralized models. The results indicate that the centralized model ($P$) provides solutions which have lower costs compared to the decentralized model ($Q$). The  gap presented varies between 2 and 6\%. These results point to the estimation errors  when assuming centralized decision making. The decentralized model, which is more realistic, leads to higher supply chain costs. 

Based on the results of the sensitivity analysis, both models suggest that biomass blends should consist of 30 to 45\% SN residues and 0\% MSW. However, the suggested percentages for the rest of biomass types differ. This result  indicates that  approximating decentralized supply chain models with centralized models leads to errors in cost estimation and less than optimal blending strategies.   
  
 \section{Summary and Conclusions}
This paper proposes models which identify blends  of biomass materials with different physical or chemical properties to satisfy the requirements of the Thermochemical conversion platform at the minimum cost. We consider requirements such as, ash  and thermal  contents. Meeting these process requirements all the time is a challenge because ash and thermal contents of biomass are random and vary by supplier. Nevertheless, these requirements should be met most  (i.e., 80-90\%) of the time in order to optimize the performance of the conversion process. We model these process requirements using chance constraints.

We consider two problem settings, centralized and decentralized  supply chains.  Most of the  supply chain literature assumes centralized system where a single decision maker has full control. However, this is not typically the case in many supply chains. We model a decentralized supply chain where the biorefinery and suppliers are independent entities  who have their own goals and objectives. We propose a  Stackelberg game where the biorefinery is the leader of this game and suppliers are the followers. We model this game using a bilevel optimization model with chance constraints.

We use the SAA  to approximate the chance constraints. We propose an exact approach and a heuristic to solve the bilevel optimization model. We also develop a solution algorithm to solve the mixed-integer programming formulation of the centralized model. We test the performance of these algorithms using a case study developed with data from South Carolina.   

Our numerical analysis indicates that the proposed SAA Algorithm finds high quality solutions in a reasonable amount of time.  The results from solving the decentralized model indicate that the blends identified contain mainly SN residues and pine. MSW is not used due to its high ash content. The amount of hybrid poplar, mixed residues, C\&D waste increases with demand for biomass. 

The blends identified when solving the centralized model use mainly pine, SP and SN residues. These blends do not contain MSW waste, and the amount of hybrid poplar, mixed residues and C\&D waste used increases with demand for biomass. Comparing the costs of meeting demand in the centralized and decentralized supply chains, one can observe that costs of the centralized model is between 2 and 6\% lower. This result indicates that, assuming centralized setting, leads to underestimates of supply chain costs. 
\newpage 
\section*{APPENDIX A}\label{APP_A}
\begin{algorithm}[!htp]
\scriptsize
\caption{SAA Algorithm}
{\bf Notation:}  $\lambda^l$ and $\lambda^u$ are the lower and upper bounds of $\lambda$; $\mu^l$ and $\mu^u$ are the lower and upper bounds of t $\mu$;  $\epsilon$ and $\delta$ be small positive constants. 
\label{Algo:binary}
\begin{algorithmic}[1]
\While{true}
\State Set $\lambda \gets \frac{\lambda^l + \lambda^u}{2}$, $\mu \gets \frac{\mu^l + \mu^u}{2}$.
\State Solve model ($\bar{\bar{P}}$) [Solve model ($\bar{Q}$)] using {\bf GUROBI} [using {\bf Heuristic Algorithm}].
\State Let $\widehat{\mathcal{W}}_s$ and $\widehat{\mathcal{J}}_s$ be the incumbent solution of $\mathcal{W}_s$ and $\mathcal{J}_s$, respectively.
\State Set $C_1 \gets 0$, $C_2 \gets 0$
\For {$s \in [1,\ldots, N]$}
\If{$\widehat{\mathcal{W}}_s > 0$}
\State $C_1 \gets C_1+1$
\EndIf
\If{$\widehat{\mathcal{J}}_s > 0$}
\State $C_2 \gets C_2+1$
\EndIf
\EndFor
\If{$C_1\geq\hat{\beta} N + \epsilon$}
\State $\lambda^l \gets \frac{\lambda^l+\lambda^u}{2}$
\ElsIf{$C_1\leq\hat{\beta} N - \epsilon$}
\State $\lambda^u \gets \frac{\lambda^l+\lambda^u}{2}$
\EndIf
\If{$C_2\geq\hat{\gamma} N+\epsilon$}
\State $\mu^l \gets \frac{\mu^l+\mu^u}{2}$
\ElsIf{$C_2\leq\hat{\gamma} N-\epsilon$}
\State $\mu^u \gets \frac{\mu^l+\mu^u}{2}$
\EndIf
\If{$|\lambda - \frac{\lambda^l+\lambda^u}{2}|\leq \delta \mbox{ AND } |\mu - \frac{\mu^l+\mu^u}{2}|\leq \delta$ }
\State break.
\EndIf
\EndWhile
\State Return $\lambda$ and $\mu$ and  solution of ($\bar{\bar{P}}$) [solution of ($\bar{Q}$)].
\end{algorithmic}
\end{algorithm}%

\begin{algorithm}[!htp]
	\scriptsize
	\begin{algorithmic}
			\State	Step 0: Initialize $p = 1$; let $\overline{\mathcal{I}}_b :=  {I}, \forall b \in \mathcal{B}$; $X_{ibp}^*$ = $\emptyset$; $\mathcal{C}_{b}^*$ = $\emptyset$
			\State \hspace{0.95cm} $Z^* = \infty$: counter = 0
			\State
			\State Step 1: \textbf{If} $\overline{\mathcal{I}}_b = \emptyset, \, \forall b \in \mathcal{B}$
			\State \hspace{1.5cm} Let $\overline{\mathcal{I}}_b := {I}, \, \forall b \in \mathcal{B}$; Let $p = p+1$
			\State \hspace{0.95cm} \textbf{End} 
			\State  \hspace{0.95cm} Find $\overline{pr}_{b}$ = $\min_{i \in \overline{\mathcal{I}}_b}\{t_{ib} + \overline{c}_{pb}\}$
			\State  \hspace{0.95cm} Let $i^*_b =  \argmin_{i \in \overline{\mathcal{I}}_b}\{t_{ib} + \overline{c}_{pb}\}$ 
			\State
			\State Step 2: Set $\hat{c}_{b} = \overline{pr}_{b}$
			\State  \hspace{0.95cm}\textbf{For} $i  \in {I}$
			\State \hspace{1.5cm}\textbf{If} $\hat{c}_{b} \geq ( t_{ib} + \overline{c}_{pb})$ 
			\State \hspace{1.95cm} $\bar{X}_{ibp} = \overline{k}_{bp}$
			\State \hspace{1.5cm} \textbf{Else} $\bar{X}_{ibp} = 0$
			\State \hspace{0.95cm}\textbf{End}  
			\State
			\State Step 3:  Solve ($\tilde{Q}^L$)
			\State
			\State Step 4: \textbf{If} ($\tilde{Q}^L$)  is not feasible:
			\State \hspace{1.5cm} Let $\overline{\mathcal{I}}_b$ = $\overline{\mathcal{I}}_b\setminus i^*_b$ {\bf GoTo} Step 1
			\State \hspace{.95cm} \textbf{Else If} $\tilde{Z}^L  < Z^*$
			\State \hspace{1.5cm} $Z^* = \tilde{Z}^L$; $X_{ibp}^*$ = $\bar{X}_{ibp}$; $\mathcal{C}_{b}^* =\hat{c}_{b}$ 
			\State \hspace{1.5cm} Let counter = 1; {\bf GoTo} Step 1
			\State \hspace{.95cm} \textbf{Else If} counter != $\nu$
			\State \hspace{1.5cm} counter = counter + 1
			\State \hspace{1.5cm} Let $\overline{\mathcal{I}}_b$ = $\overline{\mathcal{I}}_b\setminus i^*_b$; {\bf GoTo} Step 1
			\State \hspace{.95cm} \textbf{Else}
			\State \hspace{1.5cm} STOP!
			\State
			\State {\bf Return} $X_{ibp}^*$ and $\mathcal{C}_{b}^*$			
	\end{algorithmic}
	\scriptsize{\caption{Heuristic Algorithm to Solve ($\bar{Q}$) for Fixed $\lambda$, $\mu$} }
\end{algorithm}
\newpage
\section*{APPENDIX B}\label{APP_B}
 \noindent  \emph{{\bf PROPOSITION 1}: The feasible region of (${\mathbb{P}}$)  is convex.}
 
{\bf Proof}: For a given $i$ and $b$,  let $\mathbb{F}_{ibp}(X_{ib}) = \lambda_{ibp} + c_{ibp}(X_{ib}-\underline{k}_{ibp})$ for $X_{ib} \in [0, \overline{k}_{ib\mathcal{P}}]$ and for all $p=1,\ldots, \mathcal{P}$. Let function $\tilde{\mathbb{F}}_{ib}(X_{ib}) = \max_{p\in\mathcal{P}}:{\mathbb{F}_{ibp}(X_{ib})}$ for $X_{ib} \in [0, \overline{k}_{ib\mathcal{P}}]$. That means, $\tilde{\mathbb{F}}_{ib}(X_{ib}) \geq \lambda_{ibp} + c_{ibp}(X_{ib}-\underline{k}_{ibp})$ for $X_{ib} \in [0, \overline{k}_{ib\mathcal{P}}]$.  Function $\tilde{\mathbb{F}}_{ib}(X_{ib})$ is piecewise linear convex since it is the maximum of set of linear functions.   

If $\tilde{\mathbb{F}}_{ib}(X_{ib})$ is a convex function, then, for some constant $s$, the set of $X_{ib}$'s which satisfies: $\tilde{\mathbb{F}}_{ib}(X_{ib}) \leq s$ is convex \citep{luenberger1984linear}. This implies that constraints \eqref{pc_2} define a convex set. The feasible region of formulation ($\mathbb{P}$) is the intersection of linear functions  \eqref{saa_mip_1} to \eqref{saa_mip_3} and the convex set defined by \eqref{pc_2}, thus, it is convex.  $\qed$

\noindent  \emph{{\bf PROPOSITION 2}: For each feasible solution of (${\bar{\bar{P}}}$) one can find a feasible solution of (${\mathbb{P}}$), and vice versa.}

{\bf Proof}: 
\emph{Let first show that a feasible solution of (${\bar{\bar{P}}}$)  is feasible for (${\mathbb{P}}$). }

Let $\tilde{X}_{ibp}, \tilde{Z}_{ibp}$ for all $i\in I, b\in \mathcal{B}, p\in\mathcal{P}$ be a feasible solution of (${\bar{\bar{P}}}$). We will  show that such a solution satisfies \eqref{pc_2}. This solution satisfies the rest of the constraints of ($\mathbb{P}$) since they are the same for (${\bar{\bar{P}}}$).

Let $\tilde{X}_{ib} = {\tilde{X}}_{ib\tilde{p}}$, where,  $\tilde{X}_{ib\tilde{p}} = \max_{p\in\mathcal{P}}\tilde{X}_{ibp}$ and $\tilde{p}=\argmax_{p\in\mathcal{P}}\tilde{X}_{ibp}$. Based on constraints  \eqref{sumofZ}, for each $i, b$ there is exactly one $\tilde{X}_{ibp}\geq 0$ for all $p\in\mathcal{P}$, thus, $\tilde{X}_{ib} \geq 0$ and  $\tilde{X}_{ib} \in [\underline{k}_{b\tilde{p}}, \overline{k}_{b\tilde{p}}]$. Let $\tilde{F}_{ib} = \lambda_{ib\tilde{p}} + c_{ib\tilde{p}}(\tilde{X}_{ib} - \underline{k}_{b\tilde{p}})$. If $\tilde{X}_{ib} = 0$, then  $\tilde{F}_{ib}=0$. One can easily see that $\tilde{X}_{ib}$ and $\tilde{F}_{ib}$ satisfy \eqref{pc_2}.  $\qed$

\emph{Let  show that a feasible solution of (${\mathbb{P}}$) is feasible for ($\bar{\bar{P}}$).}

Let $\tilde{F}_{ib}, \tilde{X}_{ib}$ for all $i\in I, b\in \mathcal{B}$ be a feasible solution of (${{\mathbb{P}}}$). We will  show that such a solution satisfies \eqref{resourcelimits_s_32} and \eqref{sumofZ}. Recall that constraints \eqref{supply_cons} are redundant.

If $\tilde{X}_{ib} >0$ for some $i \in I, b\in \mathcal{B}$, then, let $\tilde{X}_{ib} \in [\underline{k}_{b\tilde{p}}, \overline{k}_{b\tilde{p}}]$.   Set $X_{ib\tilde{p}}=\tilde{X}_{ib}$, $Z_{ib\tilde{p}}=1$ and $X_{ibp}= Z_{ibp} = 0$, for $p\in\mathcal{P}\setminus \tilde{p}$. If $\tilde{X}_{ib} =0$ for some $i \in I, b\in \mathcal{B}$, then, let $X_{ibp}= Z_{ibp} = 0$, for $p\in\mathcal{P}$. One can easily see that $\tilde{X}_{ibp}$ and $\tilde{Z}_{ibp}$ satisfy \eqref{resourcelimits_s_32} and \eqref{sumofZ}.   $\qed$

\noindent  \emph{{\bf PROPOSITION 3}: An optimal solution of (${\mathbb{P}}$) is a lower bound of ($\bar{\bar{P}}$).}

{\bf Proof}: Let ${X}_{ib}^*, {F}_{ib}^*$ for all $i \in I, b \in \mathcal{B}$ be the optimal solution of (${\mathbb{P}}$) and let $p*$ represent the  cost bracket ${X}_{ib}^*$ corresponds to. Since this is a minimization problem, ${F}_{ib}^* = \lambda_{ib{p^*}} + c_{ib{p^*}}(\tilde{X}_{ib} - \underline{k}_{b{p^*}})$.

One can use this solution to construct a feasible solution  ${X}_{ibp}^*, {Z}_{ibp}^*$ of ($\bar{\bar{P}}$) as described in Proposition \ref{prop2}. Let  $\Delta$ be the difference in objective function value of  ($\bar{\bar{P}}$) and (${\mathbb{P}}$). Thus, $\Delta = \sum_{i\in I}\sum_{b\in \mathcal{B}}(c_{ib{p^*}}\underline{k}_{ib{p^*}} - \lambda_{ib{p^*}}) >0$.

\noindent  \emph{{\bf PROPOSITION 4}: The linear relaxation of $(\bar{Q}^F_{ib})$ provides an exact solution.}

 {\bf Proof:}  
 We prove this by contradiction. Let  $(\tilde{Q}^F_{ib})$ be the linear programming relaxation of $(\bar{Q}^F_{ib})$ obtained by relaxing \eqref{nonnegative_F1b} as follows:
 \[ Z_{p} \leq 1, \hspace{0.2in} \forall p \in \mathcal{P}.\] 
  
Let assume that an optimal solution to $(\tilde{Q}^F_{ib})$ violates constraints \eqref{nonnegative_F1b}. Let assume this solution has $X_{p^1}, X_{p^2}, Z_{p^1}, Z_{p^2} >0$.  This solution also satisfies $Z_{p^1}  + Z_{p^2} \leq1$. Let assume that $X_{p^1} \leq X_{p^2}$, thus, $\underline{k}_1 \leq  \overline{k}_1 \leq \underline{k}_2 \leq  \overline{k}_2.$ The corresponding objective function value is $\tilde{Z}^F_{ib}=\overline{c} X_{p^1}  + \overline{c} X_{p^2}$. Based on constraints  \eqref{cons_bounds}, $X_{p^1} = \overline{k}_{p^1}Z_{p^1}$ and $X_{p^2} = \overline{k}_{p^2}Z_{p^2}$. Thus, $\tilde{Z}^F_{ib}=\overline{c} \overline{k}_{p^1}Z_{p^1}  + \overline{c} \overline{k}_{p^2}Z_{p^2}.$ This means, $\tilde{Z}^F_{ib}$ is a convex combination of $\overline{c} \overline{k}_{p^1}$ and $\overline{c}\overline{k}_{p^2}.$ Consider the following cases:
\begin{itemize}
\item[(a)] If $\overline{c} \overline{k}_{p^1} > \overline{c} \overline{k}_{p^2}$: the value of $\tilde{Z}^F_{ib}$ can  increase by letting $Z_{p^1} =1$ and $Z_{p^2}=0$. 
\item[(b)] If $\overline{c} \overline{k}_{p^1} < \overline{c} \overline{k}_{p^2},$ the value of $\tilde{Z}^F_{ib}$ can  increase by letting $Z_{p^1} = 0$ and $Z_{p^2}=1$. 
\item[(c)] If $\overline{c} \overline{k}_{p^1} = \overline{c} \overline{k}_{p^2},$ the value of $\tilde{Z}^F_{ib}$ remains the same by letting $Z_{p^1} =1$ and $Z_{p^2}=0$; or  $Z_{p^1} = 0$ and $Z_{p^2}=1$. 
 \end{itemize}

In  cases (a) and (b),  we can improve the objective function value of $(\tilde{Q}^F_{ib})$. This contradicts our initial assumption that the current solution is optimal.   In case (c), a solution which does not violate \eqref{nonnegative_F1b} returns the same objective function value. This proves that the linear relaxation $(\tilde{Q}^F_{ib})$ returns a  solutions which is optimal for $(\bar{Q}^F_{ib})$.  $\qed$

\noindent  \emph{{\bf LEMMA 1}: In an optimal solution to problem $(\bar{Q}^F_{ib})$, at most one $X_{p}>0$ for $p\in\mathcal{P}$. Let $p^*$ be cost bracket for which $X_{p}>0$, then,  $X_{p^*}=\overline{k}_{p^*}$ and  $X_p = 0$ for $p \in \mathcal{P} \setminus p^*$.}
 
{\bf Proof:}   Since the objective of $(\bar{Q}^F_{ib})$ is to maximize profits, then, if $c >0$ for some $p\in\mathcal{P}$, the corresponding  $X_p$ is a candidate optimal solution. Since the objective function is linear, if $c > 0$, then, $X_p = \overline{k}_p$ maximizes profits. Based on \eqref{cons_Z1}, in an optimal solution a single cost bracket is selected. Thus, if $\max_{p\in \mathcal{P}}c\overline{k}_p > 0$, then, in an optimal solution exactly one $X_p>0$ and the remaining are 0. If $\max_{p\in \mathcal{P}}c\overline{k}_p < 0$, then, in an optimal solution  $X_p = 0, \forall p\in\mathcal{P}$. If $\max_{p\in \mathcal{P}} c_p = 0$, then, let  $p^* = \argmax_{p\in \mathcal{P}} c\overline{k}_p$. Solutions with  $X_{p^*} = 0,$  $X_{p^*} = \overline{k}$, or $X_{p^*} = \underline{k}$ and $X_p = 0$ for $p \in \mathcal{P} \setminus p^*$ lead to the same objective function value of zero.  $\qed$

\section*{APPENDIX C}\label{APP_C}
{\bf Model formulations:} Model $(Q)$. 
\begin{subequations} \label{eq:min Game_DC2}
\begin{align}
({Q}): \hspace{0.6in}&  \min: Z^L = \sum_{i \in I}\sum_{b \in \mathcal{B}} \sum_{p \in \mathcal{P}} (\mathcal{C}_{b} + f_{b})X_{ibp} \\
& \qquad \text{s.t.} \nonumber \\
& \label{CCP-1} \qquad Pr\left(\sum_{i \in I}\sum_{b \in \mathcal{B}}\sum_{p \in \mathcal{P}}(\tilde{a}_{ib}-\alpha)X_{ibp} \leq  0\right) \geq 1 - \beta,\\
& \label{CCP-2} \qquad Pr\left(\sum_{i \in I}\sum_{b \in \mathcal{B}} \sum_{p \in \mathcal{P}} e_b\tilde{h}_{ib}X_{ibp}\geq \tau\right) \geq 1- \gamma,\\
& \qquad \label{profitLB2}  \sum_{b \in \mathcal{B}}\sum_{p \in \mathcal{P}}(\mathcal{C}_{b} - \overline{c}_{pb} - t_{ib})X_{ibp} \geq 0,  \hspace{1.9in} \forall i \in I, \\
& \qquad \label{Nonnegativity-Cbp}  \mathcal{C}_{b} \geq 0,  \hspace{3.44in} \forall p \in \mathcal{P} \\
& \qquad   \label{obj_s_4}\max: Z^{F}_{i} = \sum_{b \in \mathcal{B}} \sum_{p \in \mathcal{P}} (\mathcal{C}_{b} - \overline{c}_{bp} - t_{ib})X_{ibp}, \hspace{1.39in} \forall i \in I, \\
& \qquad \qquad \label{boundcon}\underline{k}_{bp}Z_{ibp}\leq X_{ibp} \leq  \overline{k}_{bp}Z_{ibp}, \hspace{2.01in} \forall b \in \mathcal{B},  p \in \mathcal{P},\\
& \qquad \qquad \sum_{p \in \mathcal{P}} Z_{ibp} = 1, \hspace{2.76in} \forall b \in \mathcal{B},\\
& \qquad \qquad X_{ibp} \geq 0, \hspace{3.01in}\forall b \in \mathcal{B}, p \in \mathcal{P}, \\
& \qquad \qquad  \label{nonnegativeF1a} Z_{ibp} \in \{0,1\}, \hspace{2.76in}\forall  b \in \mathcal{B}, p \in \mathcal{P}.
\end{align}
\end{subequations}

{\bf Model formulations:}  Model ($\bar{Q}$).  
\begin{subequations} \label{eq:min Game_DC}
\begin{eqnarray}
(\bar{Q}): \hspace{0.3in}&& \min: \sum_{i \in I}\sum_{b \in \mathcal{B}} \sum_{p \in \mathcal{P}} (\mathcal{C}_{b} + f_{b})X_{ibp}+ \lambda \sum_{s =1}^N\mathcal{W}_s + \mu \sum_{s =1}^N\mathcal{J}_s  \nonumber \\
\vspace{2in} \nonumber\\
&& \quad \text{s.t. } \nonumber\\
&& \quad\label{UN_ash_3a}\sum_{i\in I}\sum_{b \in \mathcal{B}}\sum_{p \in \mathcal{P}}(a_{ibs} - \alpha)X_{ibp} + \mathcal{V}_s - \mathcal{W}_s =  0, \hspace{1in} \forall s =1,\ldots, N,\\
&& \quad \label{UN_thermal_3a}\tau - \sum_{i\in I}\sum_{b \in \mathcal{B}} \sum_{p \in \mathcal{P}} e_{b}h_{ibs}X_{ibp} + \mathcal{U}_s - \mathcal{J}_s =  0, \hspace{1in} \forall  s = 1,\ldots, N,\\
&& \quad \label{profit_LB3}  \sum_{b \in \mathcal{B}}\sum_{p\in \mathcal{P}}(\mathcal{C}_{b} - \overline{c}_{bp} - t_{ib})X_{ibp} \geq 0,  \hspace{1.57in} \forall i \in I, \\
&& \quad \label{Nonnegativity-Cbp2}  \mathcal{C}_{b} \geq 0,  \hspace{3.1in} \forall b \in \mathcal{B}, p \in \mathcal{P} \\
&& \qquad  \label{obj_s_5a} \max: Z_{F}^{i}  = \sum_{b \in \mathcal{B}} \sum_{p\in\mathcal{P}} (\mathcal{C}_{b} - \overline{c}_{pb} - t_{ib})X_{ibp}, \hspace{1in} \forall i \in I, \\
&& \qquad \qquad \label{Bounds} \underline{k}_{bp}Z_{ibp} \leq X_{ibp} \leq \overline{k}_{bp}Z_{ibp}, \hspace{1.55in} \forall b \in \mathcal{B},  p \in \mathcal{P},\\
&& \qquad \qquad \label{Single_Bracket} \sum_{p\in\mathcal{P}} Z_{ibp} = 1, \hspace{2.33in} \forall b \in \mathcal{B},\\
&& \qquad \qquad \label{nonnegative_F4} X_{ibp} \geq 0, \hspace{2.55in}\forall b \in \mathcal{B}, p \in \mathcal{P}, \\
&& \qquad \qquad  \label{nonnegative_F5}  Z_{ibp} \in \{0, 1\}, \hspace{2.32in}\forall  b \in \mathcal{B}, p \in \mathcal{P}.
\end{eqnarray}
\end{subequations}

{\bf Model formulations:} KKT equations of the follower's problem. 

The followers' problem is a linear program for fixed values of $c_{bp}$. We replace \eqref{obj_s_5a}- \eqref{nonnegative_F4}  with the corresponding corresponding KKT conditions.
\begin{subequations} \label{KKT}
\begin{eqnarray}
&& \mbox{Constraints: }  \eqref{Bounds}, \eqref{Single_Bracket}, \eqref{nonnegative_F4} \nonumber\\
&& \qquad \qquad  \label{nonnegative_F6}  Z_{ibp} \leq 1, \hspace{2.39in}\,\,\, \forall i \in I,  b \in \mathcal{B}, p \in \mathcal{P},\\
&& \quad \label{dual_1_DE}{ (\mathcal{C}_{b} - \overline{c}_{bp} - t_{ib})+ u_{ib} - v_{ibp} + w_{ibp}  - l_{ibp} = 0  \hspace{0.7in} \,\,\, \forall i \in I, b \in \mathcal{B}, p \in \mathcal{P},} \\
&& \quad \label{dual_2_DE} \underline{k}_{bp}v_{ibp} - \overline{k}_{bp} w_{ibp} -m_{ibp} +k_{ibp} = 0 \hspace{1.35in} \forall i \in I, b \in \mathcal{B}, p \in \mathcal{P}, \\
&& \quad \label{dual_4_DE}  (-X_{ibp} + \underline{k}_{bp}Z_{ibp})v_{ibp} = 0 \hspace{1.9in} \forall i\in I, b \in \mathcal{B}, p \in \mathcal{P}, \\
&& \quad \label{dual_5_DE}  (X_{ibp} - \overline{k}_{bp}Z_{ibp})w_{ibp} = 0 \hspace{1.98in} \forall i\in I, b \in \mathcal{B}, p \in \mathcal{P}, \\
&& \quad \label{dual_6_DE}  (\sum_{p\in\mathcal{P}} Z_{ibp} - 1)\gamma_{ib} = 0 \hspace{2.23in} \forall i\in I, b \in \mathcal{B},\\
&& \quad \label{dual_7_DE} X_{ibp}l_{ibp} = 0 \hspace{2.71in}\,\,\, \, \forall i \in I, b \in \mathcal{B}, p\in \mathcal{P},\\
&& \quad \label{dual_8_DE} Z_{ibp}m_{ibp} = 0 \hspace{2.64in}\,\,\, \, \forall i \in I, b \in \mathcal{B}, p\in \mathcal{P},\\
&& \quad \label{dual_9_DE} (Z_{ibp}-1)k_{ibp} = 0 \hspace{2.34in}\,\,\, \, \forall i \in I, b \in \mathcal{B}, p\in \mathcal{P},\\
&& \quad \label{nonnegative_U} u_{ib}, \gamma_{ib} \geq 0 \hspace{2.77in}\,\,\, \, \forall i \in I, b \in \mathcal{B},\\
&& \quad \label{nonnegative_D} v_{ibp}, w_{ibp}, l_{ibp}, m_{ibp}, k_{ibp} \geq 0 \hspace{1.75in}\,\,\, \, \forall i \in I, b \in \mathcal{B}, p \in \mathcal{P}.
\end{eqnarray}
\end{subequations}

Thus, the single level optimization model is the following.  
\begin{subequations} \label{eq:min Game_DC}
\begin{eqnarray}
\hspace{0.3in}&& \min: \sum_{i \in I}\sum_{b \in \mathcal{B}} \sum_{p \in \mathcal{P}} (\mathcal{C}_{b} + f_{b})X_{ibp}+ \lambda \sum_{s =1}^N\mathcal{W}_s + \mu \sum_{s =1}^N\mathcal{J}_s  \nonumber \\
&& \quad \text{s.t. } \eqref{UN_ash_3a} - \eqref{Nonnegativity-Cbp2}, \eqref{Bounds} - \eqref{nonnegative_F4}, \eqref{nonnegative_F6}  - \eqref{nonnegative_D}. \nonumber
\end{eqnarray}
\end{subequations}

\section*{APPENDIX D: Evaluating the Performance of SAA}\label{APP_D}
\subsection*{Finding a feasible solution for ($P$):} \label{Sec:Validating} 
 The SAA literature proposes two approaches to generate feasible solutions for ($P$) \citep{LA08}. The \emph{first approach} identifies the sample size $N$ \emph{prior} to solving ($\bar{P}$). For $\hat{\beta} < \beta$, $\hat{\gamma} < \gamma$, and for $N$ large enough, the feasible region of ($\bar{P}$) is a subset of the feasible region of $(P)$. Thus, a feasible solution of ($\bar{P}$) will be feasible for (${P}$) with high probability as $N\to \infty$. These \emph{a priori} estimates yield  very large sample size $N$  which impact the size of ($\bar{P}$) and  make it impracticable to solve. The \emph{second approach} uses a smaller sample size $N$ to find a solution $\bar{x}$ of  ($\bar{P}$), and then, conducts a \emph{a posteriori} check to see if $Pr(E^1(\bar{x}, \tilde{a}) \leq 0) \geq 1-\beta$, and $Pr(E^2(\bar{x}, \tilde{h}) \leq 0) \geq 1-\gamma$.  We use the second approach and provide details below. 

Let assume $\bar{x}$ is a solution returned by {\bf SAA Algorithm}. To estimate $p^1(\bar{x})= Pr\left({E}^1(\bar{x},\tilde{a}) \leq 0\right)$ and $p^2(\bar{x})= Pr\left({E}^2(\bar{x},\tilde{h}) \leq 0 \right)$, we sample  iid values of the random problem parameters $\tilde{a}, \tilde{h}$ ($a_1,\ldots,a_{N^\prime}, h_1,\ldots, h_{N^\prime}$). This new sample of size $N^\prime$ is used to calculate estimating probabilities $\hat{p}^1_{N^\prime}(\overline{x}) $ and $\hat{p}^2_{N^\prime}(\overline{x})$ using equations  \eqref{SAA_1p} and \eqref{SAA_2p}. By the law of large numbers, probability distribution of ${p}^1(\overline{x})$ can be approximated reasonably close by a normal distribution with mean $p^1(\overline{x})$ and variance $p^1(\overline{x})\left(1-p^1(\overline{x})\right)/N^\prime$; and probability distribution of ${p}^2(\overline{x})$ can be approximated reasonably close by a normal distribution with mean $p^2(\overline{x})$ and variance $p^2(\overline{x})\left(1-p^2(\overline{x})\right)/N^\prime$. Using this approximation, one can define one-sided ($1-\delta$)-confidence interval for $p^1(\overline{x})$ and $p^2(\overline{x})$ as follows \citep{NS06} 
\[p^1(\overline{x})\leq U^1_{\delta,N^\prime}(\overline{x}),\]
\[p^2(\overline{x})\leq U^2_{\delta,N^\prime}(\overline{x}),\]
 
\noindent where $U^1_{\delta,N^\prime}(\overline{x})= \hat{p}^1_{N^\prime}(\overline{x}) + \Phi^{-1}(1-\delta)\sqrt{\hat{p}^1_{N^\prime}(\overline{x})\left(1-\hat{p}^1_{N^\prime}(\overline{x})\right)/N^\prime}$, $U^2_{\delta,N^\prime}(\overline{x})= \hat{p}^2_{N^\prime}(\overline{x}) + \Phi^{-1}(1-\delta)\sqrt{\hat{p}^2_{N^\prime}(\overline{x})\left(1-\hat{p}^2_{N^\prime}(\overline{x})\right)/N^\prime}$  and $\Phi^{-1}()$ represents the inverse cdf of standard normal distribution. Finally, in order to check the violation of each chance constraint of $(P)$ we  compare the values of $U^1_{\beta,N^\prime}(\overline{x})$ with $\beta$ and $U^2_{\beta,N^\prime}(\overline{x})$ with $\gamma$.  A solution $\overline{x}$ returned by {\bf SAA Algorithm} is feasible to the true problem ($P$) at ($1-\delta$) confidence level if we have $U^1_{\beta,N^\prime}(\overline{x},a) \leq \beta$ and $U^2_{\beta,N^\prime}(\overline{x}, h) \leq \gamma$. \label{lemma_Q}  



\subsection*{Finding a lower bound for ($P$):} \label{sec:LB} 
The SAA can also be used to compute lower bounds for $(P)$ with high confidence.  We use the approach proposed by \cite{NS06} to calculate lower bounds. We describe this approach below.

Let $N$ be the total number of observations in a sample, and let $M$ be the total number of samples generated. The following is the procedure developed to generate a lower bound. 

First, we select values for $N^1, M^1$.  Calculate:

$$ \mathbf{B}\left(\lfloor N^1\hat{\beta}\rfloor;{\beta},N^1 \right) = \sum_{i=0}^{\lfloor N^1\hat{\beta}\rfloor} \binom{N^1}{i}\beta^i(1-\beta)^{N^1-i},$$
\noindent which is the cdf of binomial distribution and represents the probability that ${E}^1({x},a_s) > 0$ in at most  $\lfloor N^1\hat{\beta}\rfloor$  of the $N^1$ observations made ($s=1,\ldots,N^1$). Let $\pi_{N^1} = \mathbf{B}\left(\lfloor N^1\hat{\beta}\rfloor;{\beta},N^1 \right).$

Second, we choose $T^1$ to be the largest number which satisfies the following 
\[ \mathbf{B}\left(T^1-1;\pi_{N^1},M^1\right) \leq \delta. \]

Here, $\delta$ is the probability that we observe at most $T^1-1$ successes from the total of $M^1$ samples. The probability of success is $\pi_{N^1}$, and a success is a sample in which at least $\lfloor N^1\hat{\beta}\rfloor$ of the $N^1$ observations made ($s=1,\ldots,N^1$) satisfy ${E}^1({x},a^s) > 0$. We follow a similar procedure to identify $M^2$ and $N^2$.

Third, we generate $M (= \max(M^1, M^2))$ independent samples; and each sample contains $N (= \max(N^1, N^2))$ observations of the random vectors $\tilde{a}_{ib}, \tilde{h}_{ib}$. For each sample, we solve problem ($\bar{P}$) to obtain the corresponding optimal objective function values $\vartheta_{N,m} (m=1,\ldots,M$). These values represent iid realizations of the random variable $\vartheta^\ast$. In order to find lower bounds for $\vartheta^\ast$, we sort these  values in a non-decreasing order, i.e. $\vartheta_{N,(1)}\leq \vartheta_{N,(2)} \leq \ldots \leq \vartheta_{N,(M)}$. It is shown that, the quantity $\vartheta_{N,(T)}$ is a lower bound to $\vartheta^\ast$ with probability at least ($1-\delta$).   



\subsection*{Evaluating the Performance of  SAA Algorithm to Solve (${{P}}$):}\label{Sec_LowerBound_Centralized}
 See  Appendix D for details of the procedure we use to generate lower and upper bounds for ($P$) via SAA, using a method  developed by \citep{LA08}. 

To evaluate the effectiveness of the SAA  in finding feasible solutions for ($P$),  we vary the risk level $\beta, \gamma$ and sample size $N$.   Tables \ref{tab:Centralize-Feasible-00-03} and \ref{tab:Centralize-Feasible-00-02} summarize the results for $\hat{\beta} = \hat{\gamma} = 0$ and  $M=10$ replications. In Table  \ref{tab:Centralize-Feasible-00-03},  $\beta = \gamma = 0.3$, and in Table \ref{tab:Centralize-Feasible-00-02},  $\beta = \gamma = 0.2$.  For each combination of $\beta, \gamma, N$, we calculate the risk of the generated solution and the cost of the feasible solutions found, i.e., those solutions which have risk less than 0.3 and 0.2 correspondingly. For a given solution $x^*$, the risk is $Pr\{E^1(x^*,\tilde{a}) \ngtr 0\}$ and $Pr\{E^2(x^*,\tilde{h}) \ngtr 0\}$. 
 We also report the  the corresponding average, minimum, maximum and sample standard deviation ($\sigma$) of the risk and the average run time over the 10 replications.
We  report the number of feasible solutions found. For these feasible solutions, we report  the average, minimum, maximum and  standard deviation of costs.   
   
\begin{table}[htp!]
	\caption{Solutions Returned by \textbf{SAA Algorithm} for (${{P}}$) ($\beta =\gamma $ = 0.3; $\hat{\beta} =\hat{\gamma}$ = 0.0)}
	\label{tab:Centralize-Feasible-00-03}
	\scriptsize
	\centering	
	\begin{tabular}{ccccccccccc}
		\toprule
		& \multicolumn{4}{c}{Solution  Risk} &\# of Feasible&\multicolumn{4}{c}{Objective Value (in \$1,000)}&Avg Run Time\\ 
		\cmidrule(lr){2-5}
		\cmidrule(lr){7-10}
		N & Avg & Min & Max & $\sigma$ & Solutions& Avg & Min & Max & $\sigma$& (sec)\\
		\cmidrule(lr){1-1}
		\cmidrule(lr){2-5}
		\cmidrule(lr){6-6}
		\cmidrule(lr){7-10}
		\cmidrule(lr){11-11}		
		50&	0.443	&0.321	&0.651	&0.101	&0& -&-&-&-&0.13\\
		75&	0.323	&0.230	&0.468	&0.073	&3& 29,262& 29,234& 29,308 &40&0.12\\
		100	&0.274	&0.164	&0.366	&0.057	&7& 29,363& 29,304& 29,516 &71&0.16\\
		125	&0.233	&0.168	&0.294	&0.044	&10& 29,396&29,301&29,454&48&0.20\\
		150	&0.206	&0.161	&0.283	&0.037	&10&29,421&29,353&29,480&42&0.21\\
		175	&0.164	&0.121	&0.203	&0.030	&10&29,461&29,387&29,532&43&0.24\\
		200	&0.170	&0.106	&0.230	&0.038	&10&29,479&29,404&29,574&62&0.29\\
		250	&0.136	&0.099	&0.172	&0.026	&10&29,519&29,461&29,566&43&0.36\\
		300	&0.126	&0.086	&0.176	&0.029	&10&29,538&29,479&29,601&42&0.42\\
		\bottomrule	
	\end{tabular}
\end{table}

Based on the results from Table \ref{tab:Centralize-Feasible-00-03}, the algorithm finds 10 feasible solutions for $N \geq125$. Based on the results from Table \ref{tab:Centralize-Feasible-00-02}, the algorithm finds 10 feasible solutions for $N\geq 250$. The corresponding objective function values increase with $N$.

\begin{table}[htp!]
	\caption{Solutions Returned by \textbf{SAA Algorithm} for (${{P}}$) ($\beta =\gamma $ = 0.2; $\hat{\beta} =\hat{\gamma}$ = 0.0)}
	\label{tab:Centralize-Feasible-00-02}
	\scriptsize
	\centering	
	\begin{tabular}{ccccccccccc}
		\toprule
		& \multicolumn{4}{c}{Solution  Risk} &\# of Feasible&\multicolumn{4}{c}{Objective Value (in \$1,000)}&Avg Run Time\\ 
		\cmidrule(lr){2-5}
		\cmidrule(lr){7-10}
		N & Avg & Min & Max & $\sigma$ & Solutions& Avg & Min & Max & $\sigma$& (sec)\\
		\cmidrule(lr){1-1}
		\cmidrule(lr){2-5}
		\cmidrule(lr){6-6}
		\cmidrule(lr){7-10}		
		\cmidrule(lr){11-11}	
		
		50&	0.443	&0.321	&0.651	&0.101	&0&-&-&-&-&0.13\\
		75&	0.323	&0.230	&0.468	&0.073	&0&-&-&-&-&0.12\\
		100	&0.274	&0.164	&0.366	&0.057	&1&29,344&29,344&29,344&-&0.16\\
		125	&0.233	&0.168	&0.294	&0.044	&3&29,446&29,433&29,454&12&0.20\\
		150	&0.206	&0.161	&0.283	&0.037	&5&29,419&29,353&29,455&45&0.21\\
		175	&0.164	&0.121	&0.203	&0.030	&9&29,456&29,387&29,532&42&0.24\\
		200	&0.170	&0.106	&0.230	&0.038	&8&29,481&29,404&29,574&68&0.29\\
		250	&0.136	&0.099	&0.172	&0.026	&10&29,519&29,461&29,566&43&0.36\\
		300	&0.126	&0.086	&0.176	&0.029	&10&29,538&29,479&29,601&42&0.42\\

		\bottomrule	
	\end{tabular}
\end{table}

Tables  \ref{tab:Centralize-Feasible-03-03} and \ref{tab:Centralize-Feasible-02-02} summarize the results of SAA Algorithm for $\hat{\beta}$ = $\beta$ and $\hat{\gamma}$ = $\gamma$. In this case, in order to find a feasible solution to model ($P$), we increase the sample size $N$. The costs of solutions found are about 2\% lower (better) than the costs of solutions found when $\hat{\beta} = \hat{\gamma} =0.$ However, such a small improvement in solution quality has a great impact in increasing computation time due to larger sample size $N$.  
 

Tables \ref{tab:Centralize-LB-03-03} and \ref{tab:Centralize-LB-02-02} summarize the values of the lower bounds, error gaps and running time for different  values of sample size $N$ and $\hat{\beta}$, $\hat{\gamma}$. The error gap presented is smaller than 0.10\%. This error gap is reduced as sample size $N$ increases to $10,000$ or $20,000$, and for $\hat{\beta} = \hat{\gamma} = 0.2$ and $\hat{\beta} = \hat{\gamma} = 0.3$. Table \ref{tab:Centralize-LB-00-03-02} summarizes the best lower bound found for $\beta = 0.2$ and 0.3 and $\gamma  = 0.2$ and $0.3$. In these experiments,  $\hat{\beta} = \hat{\gamma} = 0.00$. The running time of the SAA Algorithm is much shorter as compared to solving the problems for $\hat{\beta} = \beta$ and $\hat{\gamma} =\gamma$ since the minimum sample size for which we can find a feasible solution, is smaller. The solutions found are within $2\%$ error gap. 
 
 Based on the results of this analysis, we decided to set $\hat{\beta} = \hat{\gamma} = 0.0$ and $N=125$ (for $\beta$ = $\gamma$ = 0.3), $N=250$ (for $\beta$ = $\gamma$ = 0.3) in our sensitivity analysis. Doing so, we  get solutions of high quality in a reasonable amount of time.  
  
\begin{table}[htp!]
	\caption{Solutions Returned by \textbf{SAA Algorithm} for (${{P}}$) ($\beta =\gamma  =\hat{\beta} =\hat{\gamma}  = 0.3$)}
	\label{tab:Centralize-Feasible-03-03}
	\scriptsize
	\centering	
	\begin{tabular}{ccccccccccc}
		\toprule
		 & \multicolumn{4}{c}{Solution  Risk} &  \# of Feasible &\multicolumn{4}{c}{Objective Value (in \$1,000)}&Average Run\\ 
		\cmidrule(lr){2-5}
		\cmidrule(lr){7-10}
		N & Avg & Min & Max & $\sigma$ & Solutions& Avg & Min & Max & $\sigma$&Time (sec.)\\
		\cmidrule(lr){1-1}
		\cmidrule(lr){2-5}
		\cmidrule(lr){6-6}
		\cmidrule(lr){7-10}
		\cmidrule(lr){11-11}				
		1,000&	 0.317& 	 0.307 	& 0.327 	 &0.006 	&0&-&-&-&-&45\\
		2,000&	 0.308 &	 0.299 	 &0.318 	 &0.007 	&2& 28,991 	 &28,985& 	 28,996& 8&123\\		
		3,000&	 0.305 &	 0.298 	 &0.312 	 &0.005 	&2& 28,992 	 &28,987 	& 28,997 	& 7&244\\
		5,000&	 0.306 	& 0.302 	 &0.314 	 &0.004 	&0&-&-&-&-&600\\
		7,500&	 0.304 &	 0.301 	 &0.308 	 &0.002 	&0&-&-&-&-&1,265\\
		10,000&	 0.303 	& 0.299 	 &0.308 	 &0.003 	&2&28,988&28,986&28,990&3&2,150\\
		20,000&	 0.302 &	 0.298 	 &0.304 	 &0.002 	&2& 28,986 	& 28,984 	& 28,987 	 &2&7,898\\			
		\bottomrule	
	\end{tabular}
\end{table}

\begin{table}[htp!]
	\caption{Solutions Returned by \textbf{SAA Algorithm} for (${{P}}$) ($\beta =\gamma  =\hat{\beta} =\hat{\gamma}  = 0.2$)}
	\label{tab:Centralize-Feasible-02-02}
	\scriptsize
	\centering	
	\begin{tabular}{ccccccccccc}
		\toprule
		  & \multicolumn{4}{c}{Solution  Risk} &\# of Feasible  &\multicolumn{4}{c}{Objective Value (in \$1,000)}&Average Run\\ 
		\cmidrule(lr){2-5}
		\cmidrule(lr){7-10}
		N & Avg & Min & Max & $\sigma$ & Solutions& Avg & Min & Max & $\sigma$&Time (sec.)\\
		\cmidrule(lr){1-1}
		\cmidrule(lr){2-5}
		\cmidrule(lr){6-6}
		\cmidrule(lr){7-10}	
		\cmidrule(lr){11-11}	
		1,000&	 0.217 &	 0.210 &	 0.225& 	 0.006 &	0&-&-&-&-&44\\
		2,000&	 0.209 	& 0.200 	& 0.219 &	 0.007 	&0&-&-&-&-&130\\
		3,000&	 0.206 	 &0.202 	& 0.213 &	 0.004 	&0&-&-&-&-&243\\
		5,000&	 0.206 	& 0.203 	& 0.209 &	 0.002 	&0&-&-&-&-&607\\
		7,500&	 0.203 	& 0.201 	& 0.206 &	 0.002 	&0&-&-&-&-&1,440\\
		10,000&	 0.204 	& 0.200 	& 0.208 &	 0.002 	&1&29,147&29,147&29,147&-&2,251\\
		20,000&	 0.202 	& 0.200 	& 0.205 &	 0.002 	&1&29,152&29,152&29,152&-&8,439\\
		\bottomrule	
	\end{tabular}
\end{table}

\begin{table}[htp!]
	\caption{Lower Bounds for (${{P}}$) ($\beta =\gamma  =\hat{\beta} =\hat{\gamma}  = 0.3$)}
	\label{tab:Centralize-LB-03-03}
	\scriptsize
	\centering	
	\begin{tabular}{ccccc}
		\toprule
		N & 2,000 & 3,000 & 10,000 & 20,000\\
		\cmidrule(lr){1-1}
		\cmidrule(lr){2-2}
		\cmidrule(lr){3-3}
		\cmidrule(lr){4-4}
		\cmidrule(lr){5-5}		
		LB Value & 28,956,315 &	28,966,006 &	28,977,520 &	28,978,728\\
		GAP ($\%$) & 0.10 & 0.07 & 0.03& 0.02\\
		Average Run&&&&\\
		Time (sec.)&123&244&2,150&7,898\\
		\bottomrule
	\end{tabular}
\end{table}

\begin{table}[htp!]
	\caption{Lower Bounds for (${{P}}$) ($\beta =\gamma  =\hat{\beta} =\hat{\gamma}  = 0.2$)}
	\label{tab:Centralize-LB-02-02}
	\scriptsize
	\centering	
	\begin{tabular}{ccc}
		\toprule
		N & 10,000 & 20,000\\
		\cmidrule(lr){1-1}
		\cmidrule(lr){2-2}
		\cmidrule(lr){3-3}
		LB Value &29,137,519&29,140,580\\
		GAP ($\%$) & 0.03 & 0.04\\
		Average Run&&\\
		Time (sec.)&2,251&8,439\\
		\bottomrule
	\end{tabular}
\end{table}

\begin{table}[htp!]
	\caption{Lower Bounds for (${{P}}$)}
	\label{tab:Centralize-LB-00-03-02}
	\scriptsize
	\centering	
	\begin{tabular}{ccc}
		\toprule
		$\hat{\beta}$ = 0.0 & $\beta$ = 0.3&  $\beta$ = 0.2\\
		$\hat{\gamma}$ = 0.0 & $\gamma$ = 0.3&  $\gamma$ = 0.2\\
		\cmidrule(lr){1-1}
		\cmidrule(lr){2-2}
		\cmidrule(lr){3-3}
		N & 125 & 250\\
		\midrule
		LB Value & 28,978,728 &	29,140,580\\
		GAP ($\%$) & 1.86 & 1.10\\
		Average Run&&\\
		Time (sec.)&0.20&0.36\\
		\bottomrule
	\end{tabular}
\end{table}

\subsection*{Evaluating the Performance of SAA Algorithm to Solve ($Q$):}
To evaluate the effectiveness of the SAA  in finding feasible solutions for ($Q$),  we vary the risk level $\beta, \gamma$ and sample size $N$. For each combination of $\beta, \gamma$  and $N$ we generate and solve 10 problems using the SAA. Table \ref{tab:Decentralize-Feasible-00-03} summarizes the results for instances with $\beta = \gamma = 0.3$, and Table \ref{tab:Decentralize-Feasible-00-02} summarizes the results for instances with $\beta = \gamma = 0.2$. 

\begin{table}[htp!]
	\caption{Solutions returned by \textbf{SAA Algorithm} for $\beta =  \gamma = 0.3$ ($\hat{\beta}$ = $\hat{\gamma}$ = 0.0)}
	\label{tab:Decentralize-Feasible-00-03}
	\scriptsize
	\centering	
	\begin{tabular}{ccccccccccc}
		\toprule
		&&&&&  &&&&&\\
		& \multicolumn{4}{c}{Solution  Risk} &\# of Feasible&\multicolumn{4}{c}{Objective Value (in \$1,000)}&Avg Run Time\\ 
		\cmidrule(lr){2-5}
		\cmidrule(lr){7-10}
		N & Ave & Min & Max & $\sigma$ & Solutions& Ave & Min & Max & $\sigma$& (sec)\\
		\cmidrule(lr){1-1}
		\cmidrule(lr){2-5}
		\cmidrule(lr){6-6}
		\cmidrule(lr){7-10}
		\cmidrule(lr){11-11}			
		30&0.550&0.439&0.697&0.068&0&-&	-&	-&	-&182\\
		50&0.420&0.286&0.529&0.082&1&30,370&30,370& 30,370&-&322\\
		75&0.317&0.241&0.407&0.044&3&30,389&30,316&30,452&68&472\\
		100&0.262&0.181&0.323&0.047&8& 30,447&30,351&30,571&78&554\\
		125&0.240&0.191&0.279&0.033&10& 30,463&30,402&30,505&36&708\\
		150&0.208&0.162&0.254&0.030&10& 30,496&30,440&30,531&32&895\\
		175&0.180&0.145&0.225&0.026&10& 30,517&30,473&30,579&39&1,137\\
		200&0.174&0.128&0.224&0.032&10& 30,544&30,459&30,631&58&1,321\\
		250&0.136&0.081&0.168&0.026&10& 30,585&30,498&30,646&50&1,940\\
		300&0.115&0.094&0.138&0.016&10& 30,608&30,564&30,654&30&2,742\\
		\bottomrule	
	\end{tabular}
\end{table}

Based on the results of Tables  \ref{tab:Decentralize-Feasible-00-03},  when $\hat{\beta} = \hat{\gamma} = 0.0$, we get 10 feasible solutions for $N \geq 125$. As $N$ increases, the number of constraints in ($\bar{Q}$) increases. The corresponding feasible region becomes smaller. This increases the likelihood that solutions found by solving  ($\bar{Q}$), are feasible to ($Q$) at this particular risk level.  Additionally, increasing $N$ results in an increase of the cost of feasible solutions found. Based on the results in Table \ref{tab:Decentralize-Feasible-00-02}, when  $\hat{\beta} = \hat{\gamma} = 0.0$, we get 10 feasible solutions for $N \geq250$. 

\begin{table}[htp!]
	\caption{Solutions returned by \textbf{SAA Algorithm} for $\beta = \gamma = 0.2$ ($\hat{\beta}$ = $\hat{\gamma}$ = 0.0)}
	\label{tab:Decentralize-Feasible-00-02}
	\scriptsize
	\centering	
	\begin{tabular}{ccccccccccc}
		\toprule
		 &&&&& &&&&&\\
		& \multicolumn{4}{c}{Solution  Risk} &\# of Feasible&\multicolumn{4}{c}{Objective Value (in \$1,000)}&Avg Run Time\\ 
		\cmidrule(lr){2-5}
		\cmidrule(lr){7-10}
		N & Avg & Min & Max & $\sigma$ & Solutions& Avg & Min & Max & $\sigma$& (sec)\\
		\cmidrule(lr){1-1}
		\cmidrule(lr){2-5}
		\cmidrule(lr){6-6}
		\cmidrule(lr){7-10}
		\cmidrule(lr){11-11}		
		30&0.550&0.439&0.697&0.068&0&-&	-&	-&	-&182\\
		50&0.420&0.286&0.529&0.082&0&-&	-&	-&	-&322\\
		75&0.317&0.241&0.407&0.044&0&-&	-&	-&	-&472\\
		100&0.262&0.181&0.323&0.047&2& 30,563	& 30,556 &	 30,571 	& 11&554\\
		125&0.240&0.191&0.279&0.033&1& 30,504 &	 30,504 &	 30,504  &	-&708\\
		150&0.208&0.162&0.254&0.030&5& 30,498 &	 30,460 &	 30,531 &28&895\\
		175&0.180&0.145&0.225&0.026&8& 30,526 &	 30,473 &	 30,579 &38&1,137\\
		200&0.174&0.128&0.224&0.032&7& 30,567 &	 30,501&	 30,631 &48&1,321\\
		250&0.136&0.081&0.168&0.026&10& 30,585 	& 30,498 &	 30,646&50&1,940\\
		300&0.115&0.094&0.138&0.016&10& 30,608 	& 30,564 &	 30,654&30&2,742\\
		\bottomrule	
	\end{tabular}
\end{table}

For this set of problems we did not develop lower bounds. This is mainly because, finding lower bounds required solving ($\bar{Q}$) for very large values of $N$. We conducted a few tests using $N = 10,000$ and $N = 20,000$, however, because of the large size of these problems, we could not find an optimal solution.  As a result, we cannot comment on the quality of the solutions found by the SAA Algorithm for the decentralized model. However, we provide such an analysis for the centralized problem in Section \ref{Sec_LowerBound_Centralized}.

 \bibliographystyle{ormsv080} 
\bibliography{myBibFile} 

\end{document}